\documentclass[11pt]{article}
\usepackage{amsfonts}
\usepackage{color}
\usepackage{amsmath,amssymb,comment}
\usepackage{verbatim}
\newtheorem{theo}{Theorem}[section]
\newtheorem{lem}[theo]{Lemma}

\newtheorem{prop}[theo]{Proposition}

\newtheorem{defi}[theo]{Definition}

\newcommand{\mysection}[1]{\section{#1} \setcounter{equation}{0}}
    \makeatletter
\def\@fnsymbol#1{\ensuremath{\ifcase#1\or *\or \ddagger\or
   \mathsection\or \mathparagraph\or \|\or **\or \dagger\dagger
   \or \ddagger\ddagger \else\@ctrerr\fi}}
    \makeatother
\newcommand{\proof}{{\sc Proof.} \quad}
\newcommand{\proofc}{{\sc Proof} \ }
\newcommand{\be}{\begin{equation} \label}
\newcommand{\ee}{\end{equation}}
\newcommand{\bea}{\begin{eqnarray}\label}
\newcommand{\eea}{\end{eqnarray}}
\newcommand{\bas}{\begin{eqnarray*}}
\newcommand{\eas}{\end{eqnarray*}}
\newcommand{\bit}{\begin{itemize}}
\newcommand{\eit}{\end{itemize}}
\newcommand{\qed}{\hfill$\Box$ \vskip.2cm}
\newcommand{\nn}{\nonumber}
\newcommand{\R}{\mathbb{R}}

\newcommand{\N}{\mathbb{N}}
\newcommand{\pO}{\partial\Omega}

\newcommand{\eps}{\varepsilon}

\newcommand{\io}{\int_\Omega}

\newcommand{\na}{\nabla}
\newcommand{\Del}{\Delta}
\newcommand{\del}{\delta}
\newcommand{\al}{\alpha}
\newcommand{\lam}{\lambda}
\newcommand{\sig}{\sigma}
\newcommand{\vt}{\vartheta}
\newcommand{\pa}{\partial}
\newcommand{\bom}{\overline{\Omega}}
\newcommand{\Om}{\Omega}

\newcommand{\mint}{- \hspace*{-3.3mm} \int}

\newcommand{\ov}{\overline}
\newcommand{\un}{\underline}
\newcommand{\wh}{\widehat}

\newcommand{\hs}{\hspace*}
\newcommand{\sm}{\setminus}

\newcommand{\vp}{\varphi}
\newcommand{\lbal}{\left\{ \begin{array}{l}}
\newcommand{\lball}{\left\{ \begin{array}{ll}}
\newcommand{\ear}{\end{array} \right.}

\newcommand{\abs}{\\[5pt]}

\newcommand{\tm}{T_{max}}

\newcommand{\uU}{\un{U}}
\newcommand{\oU}{\ov{U}}

\newcommand{\uW}{\un{W}}
\newcommand{\oW}{\ov{W}}
\newcommand{\oz}{\ov{z}}
\newcommand{\uz}{\un{z}}
\newcommand{\hU}{\wh{U}}
\newcommand{\hW}{\wh{W}}
\newcommand{\ko}{\kappa_1}
\newcommand{\kt}{\kappa_2}
\newcommand{\lo}{\lam_1}
\newcommand{\lt}{\lam_2}
\newcommand{\muu}{\mu^{(u)}}
\newcommand{\muw}{\mu^{(w)}}
\newcommand{\mus}{\mu_\star}
\newcommand{\muS}{\mu^\star}
\newcommand{\Ms}{M_\star}
\newcommand{\MS}{M^\star}
\newcommand{\kD}{k_D}
\newcommand{\KD}{K_D}
\newcommand{\kS}{k_S}
\newcommand{\KS}{K_S}

\newcommand{\parab}{{\mathcal{P}}}
\newcommand{\qarab}{{\mathcal{Q}}}
\newcommand{\whs}{\wh{s}}
\newcommand{\wht}{\wh{t}}
\newcommand{\als}{\al_\star}
\newcommand{\alss}{\al_{\star\star}}

\newcommand{\bets}{\beta_\star}
\newcommand{\betss}{\beta_{\star\star}}

\newcommand{\gams}{\gamma_\star}
\newcommand{\gamss}{\gamma_{\star\star}}

\newcommand{\dels}{\del_\star}
\newcommand{\delss}{\del_{\star\star}}

\newcommand{\ys}{y_\star}

\newcommand{\sst}{s_\star}
\newcommand{\ssst}{s_{\star\star}}
\newcommand{\sssst}{s_{\star\star\star}}
\newcommand{\ths}{\theta_\star}
\newcommand{\thss}{\theta_{\star\star}}
\newcommand{\Ts}{T_\star}
\newcommand{\TS}{T^\star}
\begin{document}
\enlargethispage{10mm}
\title{A switch in dimension dependence of critical blow-up exponents\\
in a Keller-Segel system involving indirect signal production}
\author
{
Youshan Tao\footnote{taoys@sjtu.edu.cn}\\
{\small School of Mathematical Sciences, CMA-Shanghai, Shanghai Jiao Tong University,}\\
{\small Shanghai 200240, P.R.~China}
 \and
Michael Winkler\footnote{michael.winkler@math.uni-paderborn.de}\\
{\small Universit\"at Paderborn, Institut f\"ur Mathematik,}\\
{\small 33098 Paderborn, Germany} }
\date{}
\maketitle
\begin{abstract}
\noindent
  In bounded $n$-dimensional domains with $n\ge 3$, this manuscript considers an
  initial-boundary problem for a quasilinear chemotaxis system with indirect attractant production,
  as arising, inter alia, in the modeling of effects due to phenotypical heterogeneity in microbial populations.
  Under the assumption that the rates $D$ and $S$ of diffusion and cross-diffusion are suitably regular functions
  of the population density, essentially exhibiting asymptotic behavior of the form
  \bas
	D(\xi) \simeq \xi^{m-1}
	\quad \mbox{and} \quad
	S(\xi) \simeq \xi^\sig,
	\qquad \xi \simeq \infty,
  \eas
  the identity
  \bas
	\sig=m-1+\frac{4}{n}
	\qquad \qquad (n\ge 3),
  \eas
  is shown to determine a critical line for the occurrence of blow-up.
  This considerably differs from low-dimensional cases, in which the relation
  \bas
	\sig=m+\frac{2}{n}
	\qquad \qquad (n\le 2)
  \eas
  is known to play a correspondingly pivotal role.\abs
\end{abstract}
\newpage
\section{Introduction}\label{intro}
Understanding 		
chemotactic migration in increasingly complex frameworks has been forming a core objective of
biomathematical PDE analysis in the past few years.
With answers to key questions for standard Keller-Segel systems meanwhile belonging to common knowledge in this field,
the exploration of possible effects due to more intricate couplings seems yet at a rather early stage.
In particular, only quite simple model classes have so far been made accessible to exhaustive
identifications of parameter constellations in which self-enhancing chemoattraction substantially
overbalances dissipative system ingredients in the sense that singularity formation can be observed.\abs
As a paradigmatic example for such a well-understood family of basic models, let us recall the quasilinear Keller-Segel systems
\be{01}
	\lball
	u_t = \na \cdot (D(u)\na u) - \na\cdot (S(u)\na v), \\[1mm]
	0 = \Del v - \mu + u,
	\qquad \mu=\mint_\Om u,
	\ear
\ee
for which the literature has revealed quite a precise dichotomy between parameter ranges facilitating singularity formation
and regimes within which such phenomena are absent.
Indeed, it is known that if here the pivotal ingredients $D$ and $S$ are reasonably regular functions with algebraic behavior
at large population densities, suitably generalizing
\be{proto}
	D(\xi) \simeq \xi^{m-1}
	\quad \mbox{and} \quad
	S(\xi) \simeq \xi^\sig,
	\qquad \xi \simeq \infty,
\ee
then for Neumann problems posed in bounded $n$-dimensional domains, the relation
\be{02}
	\sig=m-1+\frac{2}{n}
\ee
determines a critical line: In the case when $\sig<m-1+\frac{2}{n}$, associated initial-boundary value problems admit global bounded
solutions for arbitrary initial data, while if $\sig>m-1+\frac{2}{n}$, then some unbounded solutions exist;
for this simple model, even on the corresponding critical line where $\sig=m-1+\frac{2}{n}$ quite far-reaching knowledge
is available, inter alia on the occurrence of critical mass phenomena
(\cite{biler1998}, 
\cite{perthame_ejde2006},
\cite{cieslak_laurencot_CR2009},
\cite{cieslak_laurencot_DCDS2010}, 
\cite{nagai1995}, \cite{nagai2001}),
and significant parts of these results could in fact be extended so as to cover fully parabolic relatives of (\ref{01})
(\cite{herrero_velazquez}, \cite{senba_suzuki_AAA}, \cite{cieslak_laurencot},
\cite{taowin_subcrit}, \cite{win_collapse}, \cite{cieslak_stinner_JDE2012},
\cite{cieslak_stinner_JDE2015}, \cite{win_PLMS}).
Comparably comprehensive classifications could be achieved for some variants of (\ref{01}),
but each of these yet remain within the realm of two-component chemotaxis
systems which, thus, do not go too far beyond (\ref{01}) with regard to their level of intricacy
(\cite{tello_CPDE2022}, \cite{win_IUMJ},
\cite{win_NON}).\abs
The intention of the present study now consists in the design of an analytical approach capable of providing rigorous
evidence for a potential increase of mathematical subtlety that
is encountered when (\ref{01}) is further developed into a more complex, and hence in some application contexts more realistic,
but otherwise fairly close and natural extension.
To this end, in a bounded domain $\Om\subset\R^n$ with smooth boundary
we shall subsequently consider the initial-boundary value problem
\be{0}
	\left\{ \begin{array}{ll}	
	u_t = \na \cdot (D(u)\na u) - \na \cdot (S(u)\na v) - \ko u + \kt w,
	\qquad & x\in\Om, \ t>0, \\[1mm]
	0 = \Del v - \muw(t) + w,
	\qquad \muw(t) = \mint_\Om w,
	\qquad & x\in\Om, \ t>0, \\[1mm]
	w_t = \Del w - \lo w + \lt u
	\qquad & x\in\Om, \ t>0, \\[1mm]
	\frac{\pa u}{\pa\nu}=\frac{\pa v}{\pa\nu}=\frac{\pa w}{\pa\nu}=0,
	\qquad & x\in\pO, \ t>0, \\[1mm]
	u(x,0)=u_0(x), \quad w(x,0)=w_0(x),
	\qquad & x\in\Om,
	\end{array} \right.
\ee
which arises, as a parabolic-elliptic simplification obtained on taking fast signal diffusion limits in a standard manner
(\cite{JL}), in the modeling of chemotaxis processes that incorporate an indirect 
mechanism of attractant production, mediated through an independent diffusible component.
Recent refined studies in experimental and modeling-related literature have found circuitous taxis mechanisms of this type
to be of crucial relevance for the understanding of phenotypical heterogeneity in microbial populations
(\cite{salek}, \cite{lorenzi}, \cite{painter_JTB}, \cite{macfarlane});
we underline already here that the inclusion of the constants $\ko\ge 0$ and $\kt\ge 0$ throughout this manuscript
is mainly due to the ambition to capture all essential aspects of the setting addressed in \cite{macfarlane},
where these parameters represent rates of transition between phenotypes concentrating either on chemotaxis, or on
signal secretion.
Beyond this, systems of the form (\ref{0}) appear as special cases of attraction-repulsion models in certain
parameter ranges, with aggregation of microglia observed in Alzheimer's disease constituting a further particular
application area in this framework (\cite{fujie_senba}, \cite{tao_wang_M3AS}, \cite{ding_wang}).
Moreover, the system (\ref{0}) can be used to model the process of mass spawning in corals (\cite{coll}),
in which egg-sperm bundles release eggs, which in turn produce molecules that attract sperms in the bundles.
Independently of this, in \cite[Lemma 2]{shi_jmb2021} and \cite[Proposition 2.4]{shi_zamp2021} it is demonstrated that
some models of the form (\ref{0}) are also equivalent to certain reaction-cross-diffusion systems with
a mechanism termed "distributed delayed diffusion" there,
describing biased movement of animals due to spatial memory.\abs
As has already been observed in \cite{fujie_senba} (see also
\cite{tello_wrzosek_m3as2016}, \cite{ahn_jde2020}, \cite{sainan_wu}, \cite{sainan_wu_JMAA},
\cite{suying_liu}, \cite{mishra_wrosek_AML} for some comparable results on related systems), the
indirectness in the signal evolution process described by (\ref{0})
goes along with a markedly enhanced dissipative character that is reflected in a tendency toward relaxation
significantly stronger than that in associated boundary value problems for (\ref{01}).
In the prototypical case when $D\equiv 1$ and $S\equiv id$,
this has become manifest in a result on global existence of bounded classical solutions for arbitrarily large initial data
whenever $n\le 3$, as obtained in \cite{fujie_senba} even for a fully parabolic variant of (\ref{0}).
This contrast to (\ref{01}) has recently been further substantiated in the context of a detailed study concerned with
the general quasilinear problem (\ref{0}) under assumptions covering the choices in (\ref{proto}).
In low-dimensional situations and under the additional technical assumption that $m$ be positive,
namely, rather than (\ref{02}) the upshifted line
\be{03}
	\sig=m+\frac{2}{n}
	\qquad \qquad (n\le 2)
\ee
has been identified as being critical with respect to the occurrence of blow-up in (\ref{0}):
While singularities emerge for some solutions when $n\le 2$ and $\sig>m+\frac{2}{n}$, no such explosions can arise
when $n\le 2$ and $\sig<m+\frac{2}{n}$ (\cite{taowin_292}).\abs
{\bf Main results.} \quad
The present manuscript now confirms a noticeable switch that this critical relation undergoes when the space dimension
is enlarged. In fact, we shall see that whenever $n\ge 3$, a correspondingly pivotal role is played by the identity
\be{04}
	\sig=m-1+\frac{4}{n}
	\qquad \qquad (n\ge 3),
\ee
where in the case when $m<1-\frac{2}{n}$, in addition to this the equation $\sig=\frac{2}{n}$ will be seen to describe
a secondary critical line.
Since $-1+\frac{4}{n}<\frac{2}{n}$ for any such $n$, the change from (\ref{03}) to (\ref{04})
can roughly be interpreted as indicating that
the gain in dissipativeness obtained upon passing over from direct to indirect signal production
is abruptly reduced when the space dimension exceeds the value $n=2$ (cf.~also \cite{ding_wang}, \cite{wei_wang},
\cite{ren_cvpde2022} and \cite{wenbin_lv}
for some partial precedent results in this regard, addressing boundedness results for some fully parabolic relatives).\abs
In view of the above discussion it is not surprising that the most crucial part of this discovery will be linked to
the detection of singularity formation in respectively supercritical parameter regimes.
Here we recall that in the low-dimensional settings addressed in \cite{taowin_292}, a proof of global nonexistence
can be based on analyzing the evolution of certain moment-like functionals resembling precedents from the literature
on two-component Keller-Segel systems (\cite{biler1998}, \cite{nagai2001}).
Forming a key ingredient of the argument from \cite{taowin_292}, however, superlinear growth of $u\mapsto S(u)$,
as guaranteed by the requirements $m>0$ and $\sig>m+\frac{2}{n}$ made therein, in general fails for values of $\sig$ near the line
determined by (\ref{04}), even when $m$ is assumed to be positive.\abs
Our reasoning will accordingly be based on an entirely different approach, which at its core relies on a certain monotone structure
inherent to (\ref{0}), as reflected in a cooperativeness property of a two-component parabolic system to which
(\ref{0}) can be reduced along radial trajectories (see Section \ref{sect_coop}).
By means of an associated comparison principle, in Section \ref{sect_BU}
the occurrence of blow-up can be asserted on the basis of an explicit
construction of subsolution pairs which become singular within finite time.
Guided by the ambition to capture key features of (\ref{0}) as comprehensively as possible,
the design of these functions will be based on appropriate time-dependent splitting of the spatial domain
into a neighborhood of the origin and a corresponding outer part (Lemma \ref{lem1} and Lemma \ref{lem2}),
where in contrast to previous related approaches concerned with
simpler two-component Keller-Segel systems, 	
a significantly less singular behavior of the solution component $w$ as compared to $u$ will suggest to further
subdivide the latter outer region into an intermediate transition annulus and a time-independent outer part in our analysis
(Section \ref{sect_outer}).
In Section \ref{sect_BU_proof}, an appropriately careful adjustment of the free parameters arising in this construction
will lead to the following first of our main results.
\begin{theo}\label{theo13}
  Let $n\ge 3, R>0$ and $\Om=B_R\equiv B_R(0)\subset\R^n$, let
  \be{kl}
	\ko\ge 0,
	\quad
	\kt \ge 0,
	\quad
	\lo>0
	\quad \mbox{and} \quad
	\lt>0,
  \ee
  and assume
  \be{DSreg}
	\mbox{$D$ and $S$ belong to $C^3([0,\infty))$ with $D>0$ and $S\ge 0$ on $[0,\infty)$ and } S(0)=0,
  \ee
  an that
  \be{Du}
	D(\xi) \le \KD \xi^{m-1}
	\qquad \mbox{for all } \xi>\xi_0
  \ee
  and
  \be{Sl}
	S(\xi) \ge \kS \xi^\sig
	\qquad \mbox{for all } \xi>\xi_0
  \ee
  with some $\xi_0>0$, $\KD>0$, $\kS>0$, $m\in\R$ and $\sig\in\R$ satisfying
  \be{sig}
	\sig>m-1+\frac{4}{n}
	\qquad \mbox{and} \qquad
	\sig>\frac{2}{n}.
  \ee
  Then for any $\TS>0, \Ms>0$ and $\MS>\Ms$, one can find
  $M^{(u)}\in C^0([0,R])$ and $M^{(w)} \in C^0([0,R])$ such that
  \be{13.01}
	\limsup_{r\searrow 0} r^{-n} M^{(u)}(r) <\infty
	\qquad \mbox{and} \qquad
	\limsup_{r\searrow 0} r^{-n} M^{(w)}(r) <\infty,
  \ee
  that
  \be{13.1}
	M^{(u)}(R) \le \Ms
	\qquad \mbox{and} \qquad
	M^{(w)}(R) \le \Ms,
  \ee
  and that whenever
  \be{ir}
	u_0\in W^{1,\infty}(\Om)
	\mbox{ and }
	w_0\in W^{1,\infty}(\Om)
	\mbox{ are radially symmetric and nonnegative}
  \ee
  with
  \be{13.11}
	\Ms \le \io u_0 dx \le \MS
	\qquad \mbox{and} \qquad
	\Ms \le \io w_0 dx \le \MS
  \ee
  as well as
  \be{13.2}
	\int_{B_r(0)} u_0 dx \ge M^{(u)}(r)
	\qquad \mbox{and} \qquad
	\int_{B_r(0)} w_0 dx \ge M^{(w)}(r)
	\qquad \mbox{for all } r\in (0,R)
  \ee
  the corresponding solution $(u,v,w)$ of (\ref{0}) blows up before time $\TS$ in the following sense:
  There exist $\tm\in (0,\TS)$ as well as uniquely determined functions
  \bas
	\lbal
	u\in C^0(\bom\times [0,\tm)) \cap C^{2,1}(\bom\times (0,\tm)), \\[1mm]
	v \in C^{2,0}(\bom\times (0,\tm))
	\qquad \mbox{and} \\[1mm]
	w\in C^0(\bom\times [0,\tm)) \cap C^{2,1}(\bom\times (0,\tm))
	\ear
  \eas
  such that $u\ge 0$ and $w\ge 0$ in $\Om\times (0,\tm)$, that $\io v(\cdot,t)dx=0$ for all $t\in (0,\tm)$,
  that (\ref{0}) is satisfied in the classical sense in $\Om\times (0,\tm)$, and that
  \be{13.3}
	\limsup_{t\nearrow\tm} \|u(\cdot,t)\|_{L^\infty(\Om)}=\infty.
  \ee
\end{theo}
Characterizing criticality of the relations in (\ref{sig})
will thereafter
be completed by the derivation of results on global solvability and boundedness, widely unconditional with respect to the
initial data.
Our analysis in this direction will be arranged in such a way that, in extension of the precedent in \cite{ding_wang}
which addresses a part of the region where $\sig<m-1+\frac{4}{n}$ for a fully parabolic relative with $\ko=\kt=0$,
here the full subcritical range $\{\sig<m-1+\frac{4}{n}\} \cup \{\sig<\frac{2}{n}\}$ can be covered.
\abs
This will be achieved on the basis of essentially straightforward $L^p$ testing procedures and appropriate interpolation
arguments in Section \ref{sect_4},
where we shall firstly see that for parameter pairs below the line appearing in (\ref{04}),
not only finite-time explosions can be ruled out, but moreover also infinite-time blow-up is precluded if the parameters
$\ko,\kt,\lo$ and $\lt$ comply with the assumption (\ref{kl2}) which is particularly fulfilled if either $\ko=0$
or $(\ko,\kt)=(\lt,\lo)$ and which thus is mild enough so as to cover any of the applications
discussed above:
\begin{prop}\label{prop23}
  Let $n\ge 3$ and $\Om\subset\R^n$ be a bounded domain with smooth boundary, assume (\ref{kl}), and suppose that $D$ and $S$
  satisfy (\ref{DSreg}) as well as
  \be{Dl}
	D(\xi) \ge \kD \xi^{m-1}
	\qquad \mbox{for all } \xi>\xi_0
  \ee
  and
  \be{Su}
	S(\xi) \le \KS \xi^\sig
	\qquad \mbox{for all } \xi>\xi_0
  \ee
  with some $\xi_0>0$, $\kD>0, \KS>0$, $m\in\R$ and $\sig\in\R$ such that
  \be{sig2}
	\sig<m-1+\frac{4}{n}.
  \ee
  Then for arbitrary
  \be{init}
	u_0\in W^{1,\infty}(\Om)
	\mbox{ and }
	w_0\in W^{1,\infty}(\Om)
	\mbox{ such that $u_0\ge 0$ and $w_0\ge 0$,}
  \ee
  one can find unique
  \be{reg}
	\lbal
	u\in C^0(\bom\times [0,\infty)) \cap C^{2,1}(\bom\times (0,\infty)), \\[1mm]
	v \in C^{2,0}(\bom\times (0,\infty))
	\qquad \mbox{and} \\[1mm]
	w\in C^0(\bom\times [0,\infty)) \cap C^{2,1}(\bom\times (0,\infty))
	\ear
  \ee
  such that $u\ge 0$ and $w\ge 0$, that $\io v(\cdot,t) dx =0$ for all $t>0$, and that $(u,v,w)$ solves (\ref{0})
  in the classical sense.
  If furthermore
  \be{kl2}
	\kt \lt \le \ko \lo,
  \ee
  then there exists $C>0$ such that
  \be{23.1}
	\|u(\cdot,t)\|_{L^\infty(\Om)}
	+ \|v(\cdot,t)\|_{L^\infty(\Om)}
	+ \|w(\cdot,t)\|_{L^\infty(\Om)}
	\le C
	\qquad \mbox{for all } t>0.
  \ee
\end{prop}
We remark here that according to a simple argument detailed in \cite[Proposition 1.4]{taowin_292},
the assumption (\ref{kl2}) indeed is necessary for the boundedness property in (\ref{23.1}):
In fact, if $\kt\lt > \ko\lo$ and (\ref{DSreg}), then whenever $(u_0,w_0)$ satisfies (\ref{init}) with $u_0+w_0\not\equiv 0$,
then (\ref{0}) does not admit any global bounded classical solution.\abs
A suitable modification of the variational argument underlying Proposition \ref{prop23}
will finally
reveal that also the second inequality in (\ref{sig}) cannot be substantially weakened without changing the outcome:
If $\sig<\frac{2}{n}$, namely, then regardless of the size of $m$ a global solution can always be found.
\begin{prop}\label{prop25}
  If $n\ge 3$ and $\Om\subset\R^n$ is a bounded domain with smooth boundary, if (\ref{kl}),
  (\ref{DSreg}), (\ref{Dl}) and (\ref{Su}) hold with some
  with some $\xi_0>0$, $\kD>0, \KS>0$, $m\in\R$ and $\sig\in\R$ fulfilling
  \be{sig3}
	\sig<\frac{2}{n},
  \ee
  then whenever $u_0$ and $w_0$ comply with (\ref{init}), the problem (\ref{0}) admits a unique global classical
  solution satisfying (\ref{reg}) as well as $u(\cdot,t)\ge 0$, $w(\cdot,t)\ge 0$ and $\io v(\cdot,t)=0$ for all $t>0$.
\end{prop}
\mysection{Local existence and extensibility}
Let us begin by recording an essentially standard basic statement on local solvability and extensibility.
\begin{lem}\label{lem_loc}
  Let $n\ge 1$ and $\Om\subset\R^n$ be a bounded domain with smooth boundary,
  and assume (\ref{kl}), (\ref{DSreg}) and (\ref{init}).
  Then there exist $\tm\in (0,\infty]$ and a uniquely determined triple of functions
  \be{reg_loc}
	\lbal
	u\in C^0(\bom\times [0,\tm)) \cap C^{2,1}(\bom\times (0,\tm)), \\[1mm]
	v \in C^{2,0}(\bom\times (0,\tm))
	\qquad \mbox{and} \\[1mm]
	w\in C^0(\bom\times [0,\tm)) \cap C^{2,1}(\bom\times (0,\tm))
	\ear
  \ee
  with the properties that $u\ge 0$ and $w\ge 0$ in $\Om\times (0,\tm)$, that $\io v(\cdot,t)=0$ for all $t\in (0,\tm)$,
  that (\ref{0}) is satisfied in the classical sense in $\Om\times (0,\tm)$, and that
  \be{ext}
	\mbox{if $\tm<\infty$, \qquad then \qquad}
	\limsup_{t\nearrow\tm} \|u(\cdot,t)\|_{L^\infty(\Om)} =\infty.
  \ee
  Moreover, for each $T>0$ there exists $C(T)>0$ such that
  \be{mass_rough}
	\io u(\cdot,t) dx \le C(T)
	\qquad \mbox{for all } t\in (0,\tm) \cap (0,T),
  \ee
  and that
  \be{mass_kl2}
	\sup_{T>0} C(T) < \infty
	\qquad \mbox{if (\ref{kl2}) holds.}
  \ee
\end{lem}
\proof
  The claim concerning local existence and the extensibility criterion (\ref{ext}) can be derived by adapting
  fixed-point approaches well-established in the context of chemotaxis models
  (cf., e.g., \cite{ding_wang} for some close relative).
  The rough mass estimate (\ref{mass_rough}) readily results from an integration of the first and third equations in (\ref{0}),
  and its refinement in (\ref{mass_kl2}) can similarly be derived using the observation that the two eigenvalues of the matrix
  $A:=\left( \begin{array}{cc} - \ko & \kt \\ \lt & -\lo \end{array} \right)$ are nonpositive real numbers
  if (\ref{kl2}) holds (see \cite{taowin_292} for more details).
\qed
\mysection{Blow-up. Proof of Theorem \ref{theo13}} \label{sect_BU}
\subsection{A cooperative parabolic system satisfied by cumulated densities} \label{sect_coop}
Following classical precedents concerned with the detection of radial blow-up
in Keller-Segel type systems (\cite{JL}, \cite{biler1998}),
we will build our derivation of Theorem \ref{theo13} on an analysis of associated cumulated densities.
Unlike in the situation of two-component systems involving direct signal production, however,
in the present setting a corresponding change of variables can apparently not reduce (\ref{0}) to
a scalar parabolic problem.
In fact, a coupling between first and third solution components remains, and some basic features of this interaction
can be described as follows.
\begin{lem}\label{lem11}
  Let $n\ge 1$ and $R>0$, assume (\ref{kl}), and let $\Ms>0$ and $\MS>\Ms$.
  Then there exists $T_0=T_0^{(\Ms,\MS)}>0$ such that whenever (\ref{DSreg}) and (\ref{ir}) hold  with
  \be{11.1}
	\Ms \le \io u_0 dx \le \MS
	\qquad \mbox{and} \qquad
	\Ms \le \io w_0 dx \le \MS,
  \ee
  the corresponding solution $(u,v,w)$ of (\ref{0}) has the property that letting
  \be{U}
	U(s,t):=\int_0^{s^\frac{1}{n}} \rho^{n-1} u(\rho,t) d\rho,
	\qquad s\in [0,R^n], \ t\in [0,\tm),
  \ee
  and
  \be{W}
	W(s,t):=\int_0^{s^\frac{1}{n}} \rho^{n-1} w(\rho,t) d\rho,
	\qquad s\in [0,R^n], \ t\in [0,\tm),
  \ee
  we obtain functions $U$ and $W$ which belong to $C^1([0,R^n]\times [0,\tm)) \cap C^{2,1}((0,R^n) \times (0,\tm))$
  and satisfy $U_s\ge 0$ and $W_s\ge 0$ as well as
  \be{11.2}
	U(0,t)=W(0,t)=0,
	\quad
	U(R^n,t) \ge \frac{\mus R^n}{n}
	\quad \mbox{and} \quad
	W(R^n,t) \ge \frac{\mus R^n}{n}
	\qquad \mbox{for all $t\in (0,\tm)\cap (0,T_0)$,}
  \ee
  and for which we have
  \be{11.3}
	\parab^{(\muS)} [U,W](s,t) \ge 0
	\qquad \mbox{for all $s\in (0,R^n)$ and $t\in (0,\tm)\cap (0,T_0)$}
  \ee
  as well as
  \be{11.4}
	\qarab [U,W](s,t) =0
	\qquad \mbox{for all $s\in (0,R^n)$ and $t\in (0,\tm)\cap (0,T_0)$,}
  \ee
  where
  \be{11.5}
	\mus:=\frac{\Ms}{2|\Om|}
	\qquad \mbox{and} \qquad
	\muS:=\frac{2\MS}{|\Om|}.
  \ee
  Here and below, for $T>0$ and functions $\vp$ and $\psi$ from $C^1([0,R^n]\times [0,T))$ which satisfy $\vp_s \ge 0$ on
  $(0,R^n)\times (0,T)$ and are such that
  $\vp(\cdot,t) \in W^{2,\infty}_{loc}((0,R^n))$ and
  $\psi(\cdot,t) \in W^{2,\infty}_{loc}((0,R^n))$ for all $t\in (0,T)$, we let
  \be{P}
	\parab^{(\muS)}[\vp,\psi](s,t)
	:= \vp_t - n^2 s^{2-\frac{2}{n}} D(n\vp_s) \vp_{ss}
	- S(n\vp_s) \cdot \Big(\psi-\frac{\muS s}{n}\Big)
	+ \ko \vp
  \ee
  and
  \be{Q}
	\qarab[\vp,\psi](s,t)
	:= \psi_t - n^2 s^{2-\frac{2}{n}} \psi_{ss} + \lo \psi - \lt \vp
  \ee
  for $t\in (0,T)$ and a.e.~$s\in (0,R^n)$.
\end{lem}
\proof
  According to standard ODE theory, the linear problems
  \be{11.7}
	\lbal
	\uz_1'(t)=-\ko \uz_1(t) + \kt \uz_2(t),
	\qquad t>0, \\[1mm]
	\uz_2'(t) = -\lo \uz_2(t) + \lt \uz_1(t),
	\qquad t>0, \\[1mm]
	\uz_1(0)=\frac{\Ms}{|\Om|},
	\quad \uz_2(0)=\frac{\Ms}{|\Om|},
	\ear
  \ee
  and
  \be{11.6}
	\lbal
	\oz_1'(t)=-\ko \oz_1(t) + \kt \oz_2(t),
	\qquad t>0, \\[1mm]
	\oz_2'(t) = -\lo \oz_2(t) + \lt \oz_1(t),
	\qquad t>0, \\[1mm]
	\oz_1(0)=\frac{\MS}{|\Om|},
	\quad \oz_2(0)=\frac{\MS}{|\Om|},
	\ear
  \ee
  admit global solutions $(\uz_1,\uz_2) \in C^1([0,\infty);\R^2)$ and $(\oz_1,\oz_2)\in C^1([0,\infty);\R^2)$,
  and thanks to their continuity there exists $T_0=T_0^{(\Ms,\MS)}>0$ such that
  \be{11.8}
	\uz_1(t) \ge \frac{\Ms}{2|\Om|},
	\quad
	\uz_2(t)\ge \frac{\Ms}{2|\Om|}
	\quad \mbox{and} \quad
	\oz_2(t) \le \frac{2\MS}{|\Om|}
	\qquad \mbox{for all } t\in (0,T_0).
  \ee
  Now if (\ref{DSreg}), (\ref{ir}) and (\ref{11.1}) hold, then from (\ref{0}) it follows that
  $z_1(t):=\muu(t)\equiv \mint_\Om u(\cdot,t) dx$ and $z_2(t):=\muw(t)\equiv \mint_\Om w(\cdot,t) dx$, $t\in [0,\tm)$, satisfy
  \be{11.9}
	\lbal
	z_1'(t)=-\ko z_1(t) + \kt z_2(t),
	\qquad t\in (0,\tm), \\[1mm]
	z_2'(t) = -\lo z_2(t) + \lt z_1(t),
	\qquad t\in (0,\tm), \\[1mm]
	z_1(0)=\muu(0),
	\quad z_2(0)=\muw(0),
	\ear
  \ee
  where (\ref{11.1}) warrants that $\frac{\Ms}{|\Om|} \le z_i(0)\le \frac{\MS}{|\Om|}$ for $i\in \{1,2\}$.
  By cooperativeness of (\ref{11.9}) it thus follows from a corresponding comparison principle that
  $\uz_i(t) \le z_i(t) \le \oz_i(t)$ for all $t\in (0,\tm)$ and for $i\in \{1,2\}$, whence (\ref{11.8}) ensures that
  \be{11.10}
	\muu(t) \ge \frac{\Ms}{2|\Om|}
	\quad \mbox{and} \quad
	\frac{\Ms}{2|\Om|} \le \muw(t) \le \frac{2\MS}{|\Om|}
	\qquad \mbox{for all } t\in (0,\tm)\cap (0,T_0).
  \ee
  Since a straightforward consideration on the basis of (\ref{0}) shows that (\ref{U}) and (\ref{W}) indeed define functions
  $U$ and $W$ which enjoy the claimed regularity and monotonicity features, and which moreover satisfy
  \be{11.11}
	U(0,t)=W(0,t)=0,
	\quad
	U(R^n,t) = \frac{\muu(t) R^n}{n}
	\quad \mbox{and} \quad
	W(R^n,t) = \frac{\muw(t) R^n}{n}
	\qquad \mbox{for all } t\in (0,\tm)
  \ee
  as well as $\qarab [U,W] \equiv 0$ in $(0,R^n)\times (0,\tm)$ and
  \be{11.12}
	U_t= n^2 s^{2-\frac{2}{n}} D(nU_s) U_{ss} + S(nU_s) \cdot \Big(W-\frac{\muw(t) s}{n}\Big) - \ko U + \kt W
	\qquad \mbox{in } (0,R^n)\times (0,\tm),
  \ee
  from (\ref{11.11}) and the first two inequalities in (\ref{11.10}) we obtain (\ref{11.2}), while (\ref{11.3})
  is a consequence of (\ref{11.12}) and the last relation in (\ref{11.10}).
\qed
The cooperative character of the parabolic system associated with (\ref{11.2})-(\ref{11.4}) facilitates the derivation
of a comparison principle, the following particular version of which will form the fundament for our blow-up analysis.
\begin{lem}\label{lem_CP}
  Let $n\ge 1$, $R>0$ and $\muS>0$, assume (\ref{DSreg}) and (\ref{kl}), and suppose that for some $T>0$, the functions
  $\uU,\oU,\uW$ and $\oW$ belong to $C^1([0,R^n] \times [0,T))$ and have the properties that
  $\{\uU(\cdot,t),\oU(\cdot,t),\uW(\cdot,t),\oW(\cdot,t)\} \subset W^{2,\infty}_{loc}((0,R^n))$ for all $t\in (0,T)$,
  and that $\uU_s$ and $\oU_s$ are nonnegative on $(0,R^n)\times (0,T)$.
  If with $\parab^{(\muS)}$ and $\qarab=\qarab^{(\muS)}$ as in (\ref{P}) and (\ref{Q}) we have
  \be{C1}
	\parab^{(\muS)}[\uU,\uW](s,t) \le 0
	\quad \mbox{and} \quad
	\parab^{(\muS)}[\oU,\oW](s,t) \ge 0
	\qquad \mbox{for all $t\in (0,T)$ and a.e.~} s\in (0,R^n)
  \ee
  as well as
  \be{C2}
	\qarab[\uU,\uW](s,t) \le 0
	\quad \mbox{and} \quad
	\qarab[\oU,\oW](s,t) \ge 0
	\qquad \mbox{for all $t\in (0,T)$ and a.e.~} s\in (0,R^n),
  \ee
  and if furthermore
  \be{C3}
	\uU(0,t) \le \oU(0,t),
	\quad
	\uU(R^n,t) \le \oU(R^n,t),
	\quad
	\uW(0,t) \le \oW(0,t)
	\quad \mbox{and} \quad
	\uW(R^n,t) \le \oW(R^n,t)
  \ee
  for all $t\in [0,T)$
  as well as
  \be{C4}
	\uU(s,0) \le \oU(s,0)
	\quad \mbox{and} \quad
	\uW(s,0)\le \oW(s,0)
	\qquad \mbox{for all } s\in [0,R^n],
  \ee
  then
  \be{C5}
	\uU(s,t)\le \oU(s,t)
	\qquad \mbox{for all $s\in [0,R^n]$ and } t\in [0,T).
  \ee
\end{lem}
\proof
  We first note that for each $t\in (0,T)$, according to the inclusion
  $\{\uU(\cdot,t),\oU(\cdot,t),\uW(\cdot,t),\oW(\cdot,t)\}$ $ \subset W^{2,\infty}_{loc}((0,R^n))$
  and in line with the Rademacher theorem we can find a null set $N=N(t)\subset (0,R^n)$ such that
  $\uU_s(\cdot,t),\oU_s(\cdot,t),\uW_s(\cdot,t)$ and $\oW_s(\cdot,t)$ are differentiable in $(0,R^n)\sm N(t)$, and that
  \be{C01}
	\sup_{s\in (\eta,R^n-\eta) \sm N(t)} \Big\{ |\uU_{ss}(s,t)| + |\oU_{ss}(s,t)| + |\uW_{ss}(s,t)| + |\oW_{ss}(s,t)| \Big\}
	<\infty
	\qquad \mbox{for all } \eta \in \Big(0,\frac{R^n}{2}\Big).
  \ee
  We next fix an arbitrary $T'\in (0,T)$ and can then rely on the continuity of $\uU_s$ on $[0,R^n]\times [0,T']$ and of $S$
  on $[0,\infty)$ to pick $c_1=c_1(T')>0$ and $c_2=c_2(T')>0$ such that
  \be{C6}
	\uU_s(s,t) \le c_1
	\quad \mbox{for all } (s,t)\in [0,R^n]\times [0,T']
	\qquad \mbox{and} \qquad
	S(\xi) \le c_2
	\quad \mbox{for all } \xi\in [0,nc_1].
  \ee
  We then take any $\vt=\vt(T')>0$ such that
  \be{C7}
	\vt>\lt-\lo
	\qquad \mbox{and} \qquad
	\vt+\ko > c_2,
  \ee
  and observe that if for fixed $\eps\in (0,1)$ we let
  \be{C8}
	\vp(s,t):=\uU(s,t)-\oU(s,t) - \eps e^{\vt t},
	\qquad (s,t)\in [0,R^n]\times [0,T'),
  \ee
  then according to the continuity of $\uU$ and $\oU$ and in line with (\ref{C3}) and (\ref{C4}) we obtain that
  \be{C9}
	t_0:=\sup \Big\{ \wh{T}\in (0,T') \ \Big| \ \vp(s,t)<0 \mbox{ for all $s\in [0,R^n]$ and } t\in [0,\wh{T}]\Big\}
  \ee
  is well-defined, and that if $t_0<T'$, then there would exist $s_0\in (0,R^n)$ such that
  \be{C99}
	\vp_t(s_0,t_0)\ge 0
  \ee
  and
  \be{C10}
	0 = \vp(s_0,t_0)=\max_{s\in [0,R^n]} \vp(s,t_0).
  \ee
  According to our above selection of $N(t_0)$, the latter clearly implies that not only
  \be{C11}
	\vp_s(s_0,t_0)=0,
  \ee
  but that, apart from that, there must exist $(s_j)_{j\in\N} \subset (0,R^n)\sm N(t_0)$ such that $s_j\to s_0$ as $j\to\infty$ and
  \be{C12}
	\liminf_{j\to\infty} \vp_{ss}(s_j,t_0) \le 0
  \ee
  as well as
  \be{C13}
	\limsup_{j\to\infty} |\oU_{ss}(s_j,t_0)| < \infty.
  \ee
  Now from (\ref{C8}) and (\ref{C1}) it follows that
  \bea{C14}
	\vp_t
	&=& \uU_t - \oU_t - \vt \eps e^{\vt t} \nn\\
	&\le& n^2 s^{2-\frac{2}{n}} D(n\uU_s) \uU_{ss}
	+ S(n\uU_s) \cdot \Big( \uW - \frac{\muS s}{n}\Big) - \ko \uU \nn\\
	& & - n^2 s^{2-\frac{2}{n}} D(n\oU_s) \oU_{ss}
	- S(n\oU_s) \cdot \Big( \oW - \frac{\muS s}{n}\Big) + \ko \oU \nn\\
	& & - \vt \eps e^{\vt t} \nn\\
	&=& n^2 s^{2-\frac{2}{n}} D(n\uU_s) \vp_{ss}
	+ n^2 s^{2-\frac{2}{n}} \cdot \Big\{ D(n\uU_s) - D(n\oU_s) \Big\} \cdot \oU_{ss} \nn\\
	& & + S(n\uU_s) \cdot \Big( \uW - \frac{\muS s}{n}\Big)
	- S(n\oU_s) \cdot \Big( \oW - \frac{\muS s}{n}\Big) \nn\\
	& & - \ko \vp
	- (\vt+\ko) \eps e^{\vt t}
	\qquad \mbox{for all $t\in (0,T')$ and } s\in (0,R^n)\sm N(t),
  \eea
  where we note that by (\ref{C11}) and the continuity of $D,S,\uW$ and $\oW$,
  \bas
	D\Big(n\uU_s(s_j,t_0)\Big)
	- D\Big(n\oU_s(s_j,t_0)\Big)
	\to D\Big(n\uU_s(s_0,t_0)\Big)
	- D\Big(n\oU_s(s_0,t_0)\Big)
	=0
  \eas
  and
  \bas
	& & \hs{-20mm}
	S\Big(n\uU_s(s_j,t_0)\Big) \cdot \Big( \uW - \frac{\muS s_j}{n}\Big)
	- S\Big(n\oU_s(s_j,t_0)\Big) \cdot \Big( \oW - \frac{\muS s_j}{n}\Big) \\
	&\to&
	S\Big(n\uU_s(s_0,t_0)\Big) \cdot \Big( \uW - \frac{\muS s_0}{n}\Big)
	- S\Big(n\oU_s(s_0,t_0)\Big) \cdot \Big( \oW - \frac{\muS s_0}{n}\Big) \\
	&=& S\Big(n\uU_s(s_0,t_0)\Big) \cdot \Big( \uW(s_0,t_0)-\oW(s_0,t_0)\Big)
  \eas
  as $j\to\infty$.
  In view of (\ref{C12}), (\ref{C13}) and the continuity of $\vp_t$ and $\vp$, evaluating (\ref{C14}) at $(s_j,t_0)$ and
  taking $j\to\infty$ thus shows that
  \be{C15}
	0 \le S\Big(n\uU_s(s_0,t_0)\Big) \cdot \Big(\uW(s_0,t_0)-\oW(s_0,t_0)\Big) - (\vt+\ko)\eps e^{\vt t_0}.
  \ee
  To derive a contradiction from this, let us next make sure that
  \be{C16}
	\psi(s,t):=\uW(s,t) - \oW(s,t) - \eps e^{\vt t},
	\qquad (s,t)\in [0,R^n]\times [0,t_0],
  \ee
  satisfies
  \be{C17}
	\psi(s,t) < 0
	\qquad \mbox{for all $s\in [0,R^n]$ and } t\in [0,t_0).
  \ee
  Indeed, if this was false, then similarly to the above, as a consequence of (\ref{C3}) and (\ref{C4}) we could find
  $\wht_0\in (0,t_0), \whs_0\in (0,R^n)$ and a sequence $(\whs_j)_{j\in\N} \subset (0,R^n)\sm N(\wht_0)$ such that
  $\whs_j\to \whs_0$ as $j\to\infty$, and that
  \bas	
	\psi_t(\whs_0,\wht_0) \ge 0,
	\qquad
	\psi(\whs_0,\wht_0)=0,
	\qquad
	\psi_s(\whs_0,\wht_0)=0
	\qquad \mbox{as well as} \qquad
	\liminf_{j\to\infty} \psi_{ss}(\whs_j,\wht_0) \le 0.
  \eas
  However, since
  \bas
	\psi_t
	&=& \uW_t - \oW_t - \vt \eps e^{\vt t} \\
	&\le& n^2 s^{2-\frac{2}{n}} \uW_{ss} - \lo \uW + \lt \uU \nn\\
	& & - n^2 s^{2-\frac{2}{n}} \oW_{ss} + \lo \oW - \lt \oU
        - \vt \eps e^{\vt t} \nn\\
	&=& n^2 s^{2-\frac{2}{n}} \psi_{ss} - \lo \psi + \lt \vp
	- (\vt + \lo - \lt) \eps e^{\vt t}
	\qquad \mbox{for all $t\in (0,t_0)$ and } s\in (0,R^n) \sm N(t)
  \eas
  according to (\ref{C2}), and since $\vp(s,t)<0$ for all $(s,t)\in [0,R^n]\times [0,\wht_0]$ due to (\ref{C9}) and the
  inequality $\wht_0<t_0$, on choosing $(s,t)=(\whs_j,\wht_0)$ here and letting $j\to\infty$ this would imply that
  \bas
	0 \le \psi_t(\whs_0,\wht_0)
	\le -(\vt + \lo -\lt) \eps e^{\vt \wht_0},
  \eas
  which is incompatible with the first restriction in (\ref{C7}).\abs
  Having thereby asserted (\ref{C17}), by continuity of $\psi$ we particularly see that $\psi(s_0,t_0)\le 0$, and that hence
  (\ref{C15}) along with (\ref{C6}) entails the inequality
  \bas
	0 &\le& S\Big(n\uU_s(s_0,t_0)\Big) \cdot \eps e^{\vt t_0} - (\vt + \ko) \eps e^{\vt t_0} \\
	&\le& c_2 \cdot \eps e^{\vt t_0} - (\vt+\ko) \eps e^{\vt t_0}.
  \eas
  This contradicts the second requirement in (\ref{C7}) and thus confirms that in fact we must have $t_0=T'$, so that
  (\ref{C5}) can be concluded upon taking $\eps\searrow 0$ and then $T'\nearrow T$.
\qed
\subsection{Construction of exploding subsolutions: a general template}
In line with Lemma \ref{lem11} and Lemma \ref{lem_CP}, our strategy related to the derivation of Theorem \ref{theo13}
will at its core aim at comparison with suitable subsolution pairs $(\uU,\uW)$ for (\ref{11.2})-(\ref{11.4}) which
develop unbounded gradients at the origin within sufficiently short time intervals.
Up to a time-dependent factor to be introduced in the related final definition \ref{duw} below,
the basic structure of our candidates for this will be determined in the following two lemmata.\abs
The first of these, concentrating on the expectedly more singular component $U$, extends an increasingly steep linear
behavior near $s=0$ by an essentially power-type shape in a corresponding outer part.
We already announce here that the free parameters appearing in crucial parts of the definition (\ref{hU})
will eventually be chosen to be a suitably small positive number $\al$ and a solution $y=y(t)$ to an initiial value
problem for $y'=\gamma y^{1+\del}$ with adequately chosen $\gamma>0$ and $\del>0$ (see Lemma \ref{lem12}).
\begin{lem}\label{lem1}
  Let $n\ge 1$, $R>0$ and $\mus>0$, and let
  \be{a}
	a\equiv a^{(\mus)} := \frac{\mus R^n}{n e^\frac{1}{e} (R^n+1)}.
  \ee
  Then for all $\al\in (0,1)$ and $T>0$ and any $y\in C^1([0,T))$ with $y(t)>\frac{1}{R^n}$ for all $t\in (0,T)$, setting
  \be{hU}
	\hU(s,t) \equiv \hU^{(\mus,\al,y)} :=
	\lball
	a y^{1-\al}(t) s,
	\qquad & t\in [0,T), \ s\in \Big[ 0, \frac{1}{y(t)}\Big], \\[1mm]
	\al^{-\al} a \cdot \Big(s-\frac{1-\al}{y(t)} \Big)^\al,
	\qquad & t\in [0,T), \ s\in \Big( \frac{1}{y(t)} , R^n \Big],
	\ear
  \ee
  defines a function $\hU\in C^1([0,R^n] \times [0,T)) \cap C^0([0,T);W^{2,\infty}((0,R^n))$ which satisfies
  \be{1.1}
	\hU(0,t)=0
	\quad \mbox{and} \quad
	\hU(R^n,t) \le \frac{\mus R^n}{n}
	\qquad \mbox{for all } t\in (0,T),
  \ee
  and for which we have $\hU(\cdot,t) \in C^2([0,R^n]\sm \{\frac{1}{y(t)}\}$ for all $t\in (0,T)$, with
  \be{1.2}
	\hU_t(s,t) = \lball
	(1-\al) a y^{-\al}(t) y'(t) s,
	\qquad t\in (0,T), \ s\in (0,\frac{1}{y(t)}), \\[1mm]
	\al^{1-\al} (1-\al) a \cdot \Big(s-\frac{1-\al}{y(t)} \Big)^{\al-1} \cdot \frac{y'(t)}{y^2(t)},
	\qquad & t\in (0,T), \ s\in \Big( \frac{1}{y(t)} , R^n \Big),
	\ear
  \ee
  and
  \be{1.3}
	\hU_s(s,t) = \lball
	a y^{1-\al}(t),
	\qquad t\in (0,T), \ s\in (0,\frac{1}{y(t)}), \\[1mm]
	\al^{1-\al} a \cdot \Big(s-\frac{1-\al}{y(t)} \Big)^{\al-1},
	\qquad & t\in (0,T), \ s\in \Big( \frac{1}{y(t)} , R^n \Big),
	\ear
  \ee
  as well as
  \be{1.4}
	\hU_{ss}(s,t) = \lball
	0,
	\qquad t\in (0,T), \ s\in (0,\frac{1}{y(t)}), \\[1mm]
	- \al^{1-\al} (1-\al) a \cdot \Big(s-\frac{1-\al}{y(t)} \Big)^{\al-2},
	\qquad & t\in (0,T), \ s\in \Big( \frac{1}{y(t)} , R^n \Big).
	\ear
  \ee
\end{lem}
\proof
  The claimed regularity properties as well as the identities in (\ref{1.2})-(\ref{1.4}) can be verified in a straightforward
  fashion on the basis of (\ref{hU}).
  To confirm (\ref{1.1}), we only need to observe that according to (\ref{hU}) and the fact that $\al\ln\al \ge -\frac{1}{e}$
  and hence $\al^{-\al}= e^{-\al\ln\al} \le e^\frac{1}{e}$,
  \bas
	\hU(R^n,t)
	= \al^{-\al} a \cdot \Big(R^n-\frac{1-\al}{y(t)}\Big)^\al
	\le e^\frac{1}{e} a R^{n\al}
	\le e^\frac{1}{e} a (R^n+1)
	\qquad \mbox{for all } t\in (0,T)
  \eas
  thanks to Young's inequality and the restriction $\al<1$.
  In view of (\ref{a}), namely, (\ref{1.1}) therefore becomes evident.
\qed
With regard to their mere functional form, our candidates for the respective component related to $W$
will share the basic shape with those from Lemma \ref{lem1}.
A crucial difference, however, consists in the choice of the rate of algebraic growth in the considered outer region:
In fact, while the function $y$ in (\ref{hW}) will finally coincide with that in (\ref{hU}), the number $\beta$
will remain above the number $\frac{2}{n}$, instead of being chosen small such as $\al$ (see again Lemma \ref{lem12}).
\begin{lem}\label{lem2}
  Let $n\ge 3$, $R>0$, $\mus>0$ and $\beta\in (\frac{2}{n},1)$.
  Then whenever $T>0$ and $y\in C^1([0,T))$ is such that $y(t)>\frac{1}{R^n}$ for all $t\in (0,T)$, taking $a=a^{(\mus)}$
  from (\ref{a}) and writing
  \be{hW}
	\hW(s,t) \equiv \hW^{(\mus,\beta,y)}(s,t):= \lball
	a y^{1-\beta}(t) s,
	\qquad & t\in [0,T), \ s\in \Big[ 0, \frac{1}{y(t)}\Big], \\[1mm]
	\beta^{-\beta} a \cdot \Big(s-\frac{1-\beta}{y(t)} \Big)^\beta,
	\qquad & t\in [0,T), \ s\in \Big( \frac{1}{y(t)} , R^n \Big],
	\ear
  \ee
  we obtain an element $\hW$ of $C^1([0,R^n] \times [0,T)) \cap C^0([0,T);W^{2,\infty}((0,R^n))$ with the properties that
  \be{2.1}
	\hW(0,t)=0
	\quad \mbox{and} \quad
	\hW(R^n,t) \le \frac{\mus R^n}{n}
	\qquad \mbox{for all } t\in (0,T),
  \ee
  that
  \be{2.2}
	\hW_t(s,t) = \lball
	(1-\beta) a y^{-\beta}(t) y'(t) s,
	\qquad t\in (0,T), \ s\in (0,\frac{1}{y(t)}), \\[1mm]
	\beta^{1-\beta} (1-\beta) a \cdot \Big(s-\frac{1-\beta}{y(t)} \Big)^{\beta-1} \cdot \frac{y'(t)}{y^2(t)},
	\qquad & t\in (0,T), \ s\in \Big( \frac{1}{y(t)} , R^n \Big),
	\ear
  \ee
  that
  \be{2.3}
	\hW_s(s,t) = \lball
	a y^{1-\beta}(t),
	\qquad t\in (0,T), \ s\in (0,\frac{1}{y(t)}), \\[1mm]
	\beta^{1-\beta} a \cdot \Big(s-\frac{1-\beta}{y(t)} \Big)^{\beta-1},
	\qquad & t\in (0,T), \ s\in \Big( \frac{1}{y(t)} , R^n \Big),
	\ear
  \ee
  and that $\hW(\cdot,t) \in C^2([0,R^n]\sm \{\frac{1}{y(t)}\}$ for all $t\in (0,T)$, with
  \be{2.4}
	\hW_{ss}(s,t) = \lball
	0,
	\qquad t\in (0,T), \ s\in (0,\frac{1}{y(t)}), \\[1mm]
	- \beta^{1-\beta} (1-\beta) a \cdot \Big(s-\frac{1-\beta}{y(t)} \Big)^{\beta-2},
	\qquad & t\in (0,T), \ s\in \Big( \frac{1}{y(t)} , R^n \Big).
	\ear
  \ee
\end{lem}
\proof
  This can be seen in a way completely parallel to that from Lemma \ref{lem1}.
\qed
Upon multiplication by time-dependent factors with suitably rapid decay, introduced in order to digest ill-signed
zero-order contributions that arise in the action of $\parab^{(\muS)}$ and $\qarab$,
we hence obtain the family of candidates for subsolutions that will be used below.
\begin{defi}\label{duw}
  Let $n\ge 1$, $R>0$, $\al\in (0,1), \beta\in (\frac{2}{n},1)$, $\mus>0$ and $\theta>0$,
  and let $a=a^{(\mus)}$ be as in (\ref{a}). Given $T>0$ and
  $y\in C^1([0,T))$ such that $y(t)>\frac{1}{R^n}$ for all $t\in [0,T)$, we then let
  \be{uU}
	\uU(s,t) \equiv \uU^{(\mus,\al,\theta,y)}(s,t):=e^{-\theta t} \hU(s,t),
	\qquad s\in [0,R^n], \ t\in [0,T),
  \ee
  and
  \be{uW}
	\uW(s,t) \equiv \uW^{(\mus,\beta,\theta,y)}(s,t):=e^{-\theta t} \hW(s,t),
	\qquad s\in [0,R^n], \ t\in [0,T),
  \ee
  where $\hU=\hU^{(\mus,\al,y)}$ and
  $\hW=\hW^{(\mus,\beta,y)}$ are taken from (\ref{hU}) and (\ref{hW}).
\end{defi}
Of apparently crucial importance will be to appropriately control the effect of the
taxis-related contribution to $\parab^{(\muS)} [\uU, \uW]$ from below.
As a technical preparation for this, to be referred to in our anaylsis both near the origin and in outer regions,
the following lemma records some simple
lower bounds for the difference $\hW(s,t) - \frac{e\muS}{n} \cdot s$,
which will be used to estimate terms of the form $\hW(s,t) - \frac{e^{\theta t}\muS }{n} \cdot s$,
as
appearing in $\parab^{(\muS)} [\uU, \uW]$ (see (\ref{5.9}) and (\ref{7.11}) below).
\begin{lem}\label{lem4}
  Let $n\ge 3$, $R>0$, $\mus>0$ and $\muS>0$, and given $\beta\in (\frac{2}{n},1)$, $T>0$ and
  $y\in C^1([0,T))$ fulfilling $y(t)>\frac{1}{R^n}$ for all $t\in (0,T)$, let $a=a^{(\mus)}$ and $\hW$ be as in (\ref{a}) and
  Lemma \ref{lem2}.
  Then under the additional assumption that
  \be{4.1}
	y(t)\ge \Big(\frac{2e\muS}{na}\Big)^\frac{1}{1-\beta}
	\qquad \mbox{for all } t\in (0,T),
  \ee
  it follows that
  \be{4.2}
	\hW(s,t) - \frac{e\muS}{n} \cdot s
	\ge \frac{a}{2} y^{1-\beta}(t) s
	\qquad \mbox{for all $t\in [0,T)$ and $s\in [0,\frac{1}{y(t)}] \cap [0,\sst]$,}
  \ee
  and that
  \be{4.3}
	\hW(s,t) - \frac{e\muS}{n} \cdot s
	\ge \frac{a}{2} s^\beta
	\qquad \mbox{for all $t\in [0,T)$ and $s\in (\frac{1}{y(t)},R^n] \cap [0,\sst]$,}
  \ee
  where $\hW=\hW^{(\mus,\beta,y)}$ is as in (\ref{hW}), and where
  \be{4.4}
	\sst\equiv \sst^{(\mus,\muS,\beta)} := \bigg( \frac{n a^{(\mus)}}{2e\muS} \bigg)^\frac{1}{1-\beta}.
  \ee
\end{lem}
\proof
  Let $t\in [0,T)$. Then from (\ref{hW}) and (\ref{4.1}) we obtain that if $s\in [0,\sst]$ is such that $s\le\frac{1}{y(t)}$, then
  indeed
  \bas
	\hW(s,t) - \frac{e\muS}{n} \cdot s
	&=& \Big\{ \frac{a}{2} y^{1-\beta}(t) - \frac{e\muS}{n}\Big\}\cdot s
	+ \frac{a}{2} y^{1-\beta}(t) s \\
	&\ge& \Big\{ \frac{a}{2} \cdot \frac{2e\muS}{na} - \frac{e\muS}{n}\Big\}\cdot s
	+ \frac{a}{2} y^{1-\beta}(t) s \\
	&=& \frac{a}{2} y^{1-\beta}(t) s.
  \eas
  If $s\in [0,\sst]$ satisfies $s>\frac{1}{y(t)}$, however, then in (\ref{hW}) we can estimate
  $s -\frac{1-\beta}{y(t)}\ge s-(1-\beta)s =\beta s$, so that in this case
  \be{4.5}
	\hW(s,t) \ge \beta^{-\beta} a \cdot (\beta s)^\beta = as^\beta.
  \ee
  Since
  \bas
	\frac{a}{2} s^\beta - \frac{e\muS}{n} \cdot s
	&=& \frac{a}{2} s^\beta \cdot \Big\{ 1- \frac{2e\muS}{na} \cdot s^{1-\beta} \Big\} \\
	&\ge& \frac{a}{2} s^\beta \cdot \Big\{ 1- \frac{2e\muS}{na} \cdot \frac{na}{2e\muS} \Big\} \\[2mm]
	&=& 0
  \eas
  according to (\ref{4.4}), from (\ref{4.5}) we obtain (\ref{4.3}).
\qed
\subsection{Subsolution properties: inner region}
Now due to the simple linear structure of both $\uU$ and $\uW$ near $s=0$,
the action of both $\parab^{(\muS)}$ and $\qarab$ is not influenced by diffusion in this region.
In particular, with regard to the first of these operators our derivation of respective subsolution properties
essentially reduces to an analysis of a transport mechanism which, in the presence of suitably chosen singular functions $y$,
can be viewed as a two-component generalization of a Burgers-type shock formation.
As this mechanism needs to be driven by taxis, it seems natural that a diffusion-independent minimal
size of $\sig$ will be necessary for such an emergence of singularities.
\begin{lem}\label{lem5}
  Let $n\ge 3$ and $R>0$, and suppose that (\ref{kl}), (\ref{DSreg}) and (\ref{Sl}) hold with some
  $\xi_0>0, \kS>0$ and some
  \be{5.1}
	\sig>\frac{2}{n}.
  \ee
  Then there exist $\als \in (0,1)$, $\bets\in (\frac{2}{n},1)$  and  $\dels>0$ such that for any choice of $\mus>0$ and $\muS>0$ one can find
  $\ys=\ys(\mus,\muS)>\frac{1}{R^n}$ and $\gams=\gams(\mus,\muS)>0$ with the property that if $T>0$ and $y\in C^1([0,T))$ are
  such that
  \be{5.2}
	y(t) \ge \ys
	\qquad \mbox{for all } t\in (0,T)
  \ee
  and
  \be{5.3}
	y'(t) \le \gams y^{1+\dels}(t)
	\qquad \mbox{for all } t\in (0,T),
  \ee
  then for arbitrary $\al\in (0,\als)$ and $\beta\in (\frac{2}{n},\bets)$, and whenever $\theta>0$ is such that
  \be{5.5}
	\theta\ge \ko,
  \ee
  the functions $\uU=\uU^{(\mus,\al,\theta,y)}$ and
  $\uW=\uW^{(\mus,\beta,\theta,y)}$ from (\ref{uU}) and (\ref{uW}) satisfy
  \be{5.4}
	\parab^{(\muS)} [\uU,\uW](s,t) \le 0
	\qquad \mbox{for all $t\in (0,T)\cap (0,\frac{1}{\theta})$ and $s\in (0,\frac{1}{y(t)})$}.
  \ee
\end{lem}
\proof
  Since $(1-\al)\sig+\al-\beta \to \sig-\frac{2}{n}$ as $(\al,\beta)\to (0,\frac{2}{n})$, according to (\ref{5.1}) it is possible
  to pick $\als\in (0,\frac{1}{2}]$ and $\bets\in (\frac{2}{n},1)$ such that
  \be{5.6}
	(1-\al)\sig + \al-\beta
	\ge \dels := \frac{\sig-\frac{2}{n}}{2} >0
	\qquad \mbox{for all $\al\in (0,\als)$ and } \beta\in \Big(\frac{2}{n},\bets\Big).
  \ee
  Given $\mus>0$ and $\muS>0$, we then choose $\ys=\ys(\mus,\muS)>\frac{1}{R^n}$
  large enough such that with $a=a^{(\mus)}$ taken from
  (\ref{a}) we have
  \be{5.7}
	\ys \ge 1,
	\qquad
	\ys^\frac{1}{2} \ge \frac{e\xi_0}{na}
	\qquad \mbox{and} \qquad
	\ys \ge \Big(\frac{2e\muS}{na}\Big)^\frac{1}{1-\beta}
	\quad \mbox{for all } \beta\in \Big(\frac{2}{n},\bets\Big),
  \ee
  and we claim that then the intended conclusion holds if we let
  \be{5.8}
	\gams(\mus,\muS):=\frac{n^\sig \kS a^\sig e^{-\sig}}{2}.
  \ee
  To verify this, we fix $\al\in (0,\als), \beta\in (\frac{2}{n},\bets)$ and $\theta>0$ with $\theta\ge \ko$, and first
  go back to (\ref{P}) to see that since the function $\hU$, well-defined according to (\ref{5.2}) and the inequality
  $y>\frac{1}{R^n}$, satisfies
  $(e^{-\theta t} \hU)_t = e^{-\theta t} \hU_t - \theta e^{-\theta t} \hU$ for $t\in (0,T)$ and $s\in (0,R^n)$, it follows that
  \bea{5.9}
	& & \hs{-20mm}
	e^{\theta t} \cdot \parab^{(\muS)} [\uU,\uW](s,t) \nn\\
	&=& \hU_t - \theta \hU - n^2 s^{2-\frac{2}{n}} D\big( ne^{-\theta t} \hU_s\big) \hU_{ss}
	- S\big(ne^{-\theta t} \hU_s\big) \cdot \Big(\hW-\frac{\muS e^{\theta t} s}{n} \Big)
	+ \ko \hU \nn\\
	&\le& \hU_t
	- S\big(ne^{-\theta t} \hU_s\big) \cdot \Big(\hW-\frac{\muS e^{\theta t} s}{n} \Big)
	\qquad \mbox{for all $t\in (0,T)$ and $s\in (0,\frac{1}{y(t)})$,}
  \eea
  because $(-\theta+\ko)\hU \le 0$ according to the nonnegativity of $\hU$, and because (\ref{1.4}) asserts that
  $\hU_{ss} \equiv 0$ in the considered region.
  Here we note that the first two restrictions in (\ref{5.7}) guarantee that in line with (\ref{1.3}) and our assumption
  $\al\le\frac{1}{2}$ we have
  \bas
	n e^{-\theta t} \hU_s
	\ge \frac{n}{e} \hU_s
	= \frac{n}{e} \cdot a y^{1-\al}(t)
	\ge \frac{n}{e} \cdot a y^\frac{1}{2}(t)
	\ge \frac{n}{e} \cdot a \ys^\frac{1}{2}
	\ge \xi_0
	\qquad \mbox{for all $t\in (0,T)\cap (0,\frac{1}{\theta})$ and $s\in (0,\frac{1}{y(t)})$,}
  \eas
  so that we may rely on (\ref{Sl}) in estimating
  \bas
	S\big(ne^{-\theta t} \hU_s\big)
	&\ge& \kS \cdot \big( ne^{-\theta t} \hU_s\big)^\sig \nn\\
	&=& n^\sig \kS a^\sig e^{-\theta\sig t} y^{(1-\al)\sig}(t) \nn\\
	&\ge& n^\sig \kS a^\sig e^{-\sig} y^{(1-\al)\sig}(t)
	\qquad \mbox{for all $t\in (0,T)\cap (0,\frac{1}{\theta})$ and $s\in (0,\frac{1}{y(t)})$}
  \eas
  according to (\ref{1.3}).\abs
  Apart from that, we observe that the last inequality in (\ref{5.7}) warrants that (\ref{4.1}) holds, and that
  with $s_1=s_1^{(\beta,\mus,\muS)}$ as in (\ref{4.4}) we have
  \bas
	\frac{1}{y(t)} \le \frac{1}{\ys}
	\le  \Big(\frac{na}{2e\muS}\Big)^\frac{1}{1-\beta}
	= s_1
	\qquad \mbox{for all } t\in (0,T).
  \eas
  Therefore, Lemma \ref{lem4} applies so as to ensure that in (\ref{5.9}),
  \bas
	\hW - \frac{\muS e^{\theta t} s}{n}
	\ge \hW - \frac{\muS es}{n}
	\ge \frac{a}{2} y^{1-\beta}(t) s
	\qquad \mbox{for all $t\in (0,T)\cap (0,\frac{1}{\theta})$ and $s\in (0,\frac{1}{y(t)})$,}
  \eas
  whence using (\ref{1.2}) and the nonnegativity of $\al$ we infer from (\ref{5.9}) that
  for all $t\in (0,T)\cap (0,\frac{1}{\theta})$ and $s\in (0,\frac{1}{y(t)})$,
  \bas
	e^{\theta t} \cdot \parab^{(\muS)} [\uU,\uW](s,t)
	&\le& a y^{-\al}(t) y'(t) s
	- n^\sig \kS a^\sig e^{-\sig} y^{(1-\al)\sig}(t) \cdot \frac{a}{2} y^{1-\beta}(t) s \\\
	&=& a y^{-\al}(t) s \cdot \bigg\{ y'(t) - \frac{n^\sig \kS a^\sig e^{-\sig}}{2} \cdot y^{1+(1-\al)\sig + \al - \beta}(t)
		\bigg\}.
  \eas
  Since (\ref{5.6}) in conjunction with the inequality $y\ge 1$ guarantees that
  \bas
	y^{1+(1-\al)\sig + \al - \beta}(t) \ge y^{1+\dels}(t)
	\qquad \mbox{for all } t\in (0,T),
  \eas
  due to (\ref{5.3}) this establishes (\ref{5.4}) with $\gams$ as in (\ref{5.8}).
\qed
Yet more easily can the action of $\qarab$ on $(\uU,\uW)$ be described in this inner part:
\begin{lem}\label{lem6}
  Let $n\ge 3$ and $R>0$, assume (\ref{kl}), and suppose that with some $\al\in (0,1), \beta\in (\frac{2}{n},1)$ and $T\in (0,1)$,
  the function $y\in C^1([0,T))$ satisfies
  \be{6.1}
	y(t)>\max \Big\{ 1 \, , \, \frac{1}{R^n}\Big\}
	\quad \mbox{and} \quad
	y'(t) \le  \lt y^{1-\al+\beta}(t)
	\qquad \mbox{for all } t\in (0,T).
  \ee
  Then whenever
  \be{6.3}
	\theta\ge\lo,
  \ee
  for $\uU=\uU^{(\mus,\al,\theta,y)}$ and
  $\uW=\uW^{(\mus,\beta,\theta,y)}$ as in (\ref{uU}) and (\ref{uW}) we have
  \be{6.4}
	\qarab [\uU,\uW](s,t) \le 0
	\qquad \mbox{for all $t\in (0,T)$ and $s\in (0,\frac{1}{y(t)})$.}
  \ee
\end{lem}
\proof
  For arbitrary $\theta\ge\lo$, we recall (\ref{Q}), Lemma \ref{lem2} and Lemma \ref{lem1} to see that by nonnegativity of $\hW$,
  \bas
	e^{\theta t} \cdot \qarab [\uU,\uW](s,t)
	&=& \hW_t - \theta \hW - n^2 s^{2-\frac{2}{n}} \hW_{ss} + \lo \hW - \lt \hU \\
	&\le& \hW_t - n^2 s^{2-\frac{2}{n}} \hW_{ss} - \lt \hU \\
	&=& (1-\beta) a y^{-\beta}(t) y'(t) s
	- \lt a y^{1-\al}(t) s
	\qquad \mbox{for all $t\in (0,T)$ and $s\in (0,\frac{1}{y(t)})$.}
  \eas
  Since (\ref{a}) and (\ref{6.1}) imply that
  \bas
	\frac{(1-\beta) a y^{-\beta}(t) y'(t)}{\lt a y^{1-\al}(t)}
	= \frac{1-\beta}{\lt} \cdot y^{-1+\al-\beta}(t) y'(t)
	\le 1
	\qquad \mbox{for all } t\in (0,T),
  \eas
  this already entails (\ref{6.4}).
\qed
\subsection{Subsolution properties: outer region}\label{sect_outer}
In view of (\ref{hU}) and (\ref{hW}), in the region complementary to that in Lemma \ref{lem5} and Lemma \ref{lem6}
the influences of diffusion, now no longer trivial, need to be appropriately overbalanced by the cross-diffusive
part of $\parab^{(\muS)}[\uU,\uW]$.
To achieve this, we will separately examine the action of both parabolic operators under consideration in
a transition regime where $\frac{1}{y(t)}<s<s_0$ with some fixed and time-independent $s_0$, and an associated
neighborhood of the outer boundary where $s_0<s<R^n$.\abs
The most crucial part of our analysis in these subdomains will be contained in the following lemma which,
under appropriate assumptions on $\al,\beta$ and $y$, identifies
a number $\ssst$ with the property that $\parab^{(\muS)}[\uU,\uW] \le 0$ when $\frac{1}{y(t)} < s < \ssst$,
provided that $\sig$ lies above the threshold appearing in Theorem \ref{theo13};
the final choice of $s_0$ will accordingly be such that, inter alia, $s_0\le \ssst$.
\begin{lem}\label{lem7}
  If $n\ge 3$ and $R>0$, and if (\ref{kl}), (\ref{DSreg}), (\ref{Du}) and (\ref{Sl}) hold with some $\xi_0>0, \kD>0, \kS>0$,
  $m\in\R$ and $\sig\in\R$ such that
  \be{7.1}
	\sig>m-1+\frac{4}{n}
	\qquad \mbox{and} \qquad
	\sig>\frac{2}{n},
  \ee
  then there exist $\alss\in (0,1), \betss\in (\frac{2}{n},1)$ and $\delss>0$ with the property that whenever
  $\al\in (0,\alss), \beta\in (\frac{2}{n},\betss)$, $\mus>0$ and $\muS>0$, it is possible to fix
  $\ssst=\ssst^{(\mus,\muS,\al,\beta)} \in (0,1]$ and
  $\gamss=\gamss^{(\mus,\muS,\al,\beta)}>0$ in such a way that given any $T>0$ and $y\in C^1([0,T))$ fulfilling
  \be{7.2}
	y(t)>\max \Big\{ 1 \, , \, \frac{1}{R^n}\Big\}
	\quad \mbox{and} \quad
	0\le y'(t) \le \gamss y^{1+\delss}(t)
	\qquad \mbox{for all } t\in (0,T),
  \ee
  and any $\theta>0$ satisfying
  \be{7.5}
	\theta\ge \ko,
  \ee
  with $\uU=\uU^{(\mus,\al,\theta,y)}$ and $\uW=\uW^{(\mus,\beta,\theta,y)}$ taken from (\ref{uU}) and (\ref{uW}) we have
  \be{7.4}
	\parab^{(\muS)} [\uU,\uW](s,t) \le 0
	\qquad \mbox{for all $t\in (0,T) \cap (0,\frac{1}{\theta})$ and $s\in (\frac{1}{y(t)},R^n) \cap (0,\ssst)$.}
  \ee
\end{lem}
\proof
  Writing $\xi_+:=\max\{\xi,0\}$ for $\xi\in\R$, we note that as $(\al,\beta)\to (0,\frac{2}{n})$ we have
  \bas
	-\frac{2}{n}-\beta + (1-\al)(\sig-m+1) + \al \to -\frac{4}{n} + \sig - m +1
  \eas
  and
  \bas
	1- \Big[ \beta - (1-\al)(\sig-1)\Big]_+ \to 1- \Big[ \frac{2}{n}-\sig+1\Big]_+
	\ge \min \Big\{ 1, \sig-\frac{2}{n} \Big\},
  \eas
  and use (\ref{7.1}) to thus find $\alss\in (0,1)$, $\betss\in (\frac{2}{n},1)$ and $\delss>0$ such that
  \be{7.6}
	-\frac{2}{n}-\beta + (1-\al)(\sig-m+1) + \al >0
	\qquad \mbox{for all $\al\in (0,\alss)$ and } \beta \in \Big(\frac{2}{n},\betss\Big),
  \ee
  and that
  \be{7.7}
	1- \Big[ \beta - (1-\al)(\sig-1)\Big]_+ \ge \delss
	\qquad \mbox{for all $\al\in (0,\alss)$ and } \beta \in \Big(\frac{2}{n},\betss\Big).
  \ee
  For fixed $\al\in (0,\alss)$, $\beta\in (\frac{2}{n},\betss)$, $\mus>0$ and $\muS>0$, with
  $a=a^{(\mus)}$ and $\sst=\sst^{(\mus,\muS,\beta)}$ taken from (\ref{a}) and (\ref{4.4}) we then use (\ref{7.6}) to
  choose $\ssst=\ssst^{(\mus,\muS,\al,\beta)} \in (0,1]$ and $\gamss=\gamss^{(\mus,\muS,\al,\beta)}>0$ suitably small such that
  \be{7.8}
	\ssst \le \sst
  \ee
  and
  \be{7.9}
	\frac{n \al^{1-\al} a \ssst^{\al-1}}{e} \ge \xi_0
  \ee
  as well as
  \be{7.91}
	\frac{2 n^{m+1-\sig} \al^{(1-\al)(m-\sig)} (1-\al) a^{m-1-\sig} \KD e^{\sig-m+1}}{\kS}
		\cdot c_1(\al) \ssst^{-\frac{2}{n}-\beta+(1-\al)(\sig-m+1) + \al}
	\le \frac{1}{2}
  \ee
  and
  \be{7.92}
	\frac{2\al^{(1-\al)(1-\sig)} (1-\al) e^\sig}{n^\sig a^\sig \kS} \cdot c_2(\al) \gamss \le \frac{1}{2},
  \ee
  where
  \be{7.93}
	c_1(\al):=\max \big\{ 1 \, , \, \al^{(1-\al)(\sig-m+1)+\al-2} \big\}
  \ee
  and
  \be{7.94}
	c_2(\al):=\max \big\{ 1 \, , \, \al^{(1-\al)(\sig-1)} \big\}.
  \ee
  Now assuming that $T>0$ and $y\in C^1([0,T))$ satisfy (\ref{7.2}), and that $\theta>0$ is such that (\ref{7.5}) holds,
  we first note that in the identity
  \bea{7.11}
  \hspace*{-5mm}
	e^{\theta t} \parab^{(\muS)} [\uU,\uW](s,t)
	&=& \hU_t - \theta \hU - n^2 s^{2-\frac{2}{n}} D\big(n e^{-\theta t} \hU_s\big) \hU_{ss} \nn\\
	& & - S\big(n e^{-\theta t} \hU_s\big) \cdot \Big( \hW - \frac{\muS e^{\theta t} s}{n}\Big)
	+ \ko \hU,
	\qquad t\in (0,T), \mbox{$s\in (\frac{1}{y(t)},R^n)$,}
  \eea
  according to (\ref{7.8}) and Lemma \ref{lem4} we have
  \be{7.12}
	\hW - \frac{\muS e^{\theta t} s}{n}
	\ge \hW - \frac{\muS es}{n}
	\ge \frac{a}{2} s^\beta
	\qquad \mbox{for all $t\in (0,T)\cap (0,\frac{1}{\theta})$ and $s\in (\frac{1}{y(t)},R^n) \cap (0,\ssst)$}.
  \ee
  We moreover observe that
  \be{7.13}
	s \ge s-\frac{1-\al}{y(t)}
	\ge s-(1-\al) s = \al s
	\qquad \mbox{for all $t\in (0,T)$ and $s>\frac{1}{y(t)}$,}
  \ee
  which due to (\ref{7.9}) particularly ensures that in line with (\ref{1.3}) and the inequality $\al<1$,
  \bas
	n e^{-\theta t} \hU_s
	&\ge& \frac{n}{e} \hU_s
	= \frac{n}{e} \cdot \al^{1-\al} a \cdot \Big( s-\frac{1-\al}{y(t)} \Big)^{\al-1}
	\ge \frac{n}{e} \cdot \al^{1-\al}  a \ssst^{\al-1} \nn\\
	&\ge& \xi_0
	\qquad \mbox{for all $t\in (0,T)\cap (0,\frac{1}{\theta})$ and $s\in (\frac{1}{y(t)},R^n) \cap (0,\ssst)$}.
  \eas
  We may therefore draw on (\ref{Du}) and (\ref{Sl}) to see that thanks to (\ref{7.12}) and (\ref{7.5})
  and the fact that $\hU_{ss}\le 0$ by (\ref{1.4}), the identity in
  (\ref{7.11}) can be turned into the inequality
  \bea{7.14}
	e^{\theta t} \parab^{(\muS)} [\uU,\uW](s,t)
	&\le& \hU_t - n^2 s^{2-\frac{2}{n}} \KD \cdot \big(n e^{-\theta t} \hU_s \big)^{m-1} \hU_{ss}
	- \kS \cdot \big( n e^{-\theta t} \hU_s\big)^\sig \cdot \frac{a}{2} s^\beta \nn\\[2mm]
	& & \hs{10mm}
	\mbox{for all $t\in (0,T)\cap (0,\frac{1}{\theta})$ and $s\in (\frac{1}{y(t)},R^n) \cap (0,\ssst)$}.
  \eea
  Here by (\ref{1.3}) and (\ref{1.4}),
  \bea{7.15}
	& & \hs{-20mm}
	\frac{-n^2 s^{2-\frac{2}{n}} \KD \cdot ( n e^{-\theta t} \hU_s)^{m-1} \hU_{ss}}
		{\kS \cdot (n e^{-\theta t} \hU_s)^\sig \cdot \frac{a}{2} s^\beta} \nn\\
	&=& - \frac{2n^{m+1-\sig} \KD e^{(\sig-m+1)\theta t}}{a\kS} \cdot s^{2-\frac{2}{n}-\beta} \hU_s^{-(\sig-m+1)} \hU_{ss} \nn\\
	&=& \frac{2n^{m+1-\sig} \KD e^{(\sig-m+1)\theta t}}{a\kS} \cdot s^{2-\frac{2}{n}-\beta} \cdot
		\Big\{ \al^{1-\al} a \cdot \Big(s-\frac{1-\al}{y(t)}\Big)^{\al-1} \Big\}^{-(\sig-m+1)} \times \nn\\
	& & \hs{60mm}
		\al^{1-\al} (1-\al) a \cdot \Big(s-\frac{1-\al}{y(t)}\Big)^{\al-2} \nn\\
	&=& \frac{2n^{m+1-\sig} \al^{(1-\al)(m-\sig)} (1-\al) a^{m-1-\sig} \KD e^{(\sig-m+1)\theta t}}{\kS} \times \nn\\
	& & \hs{60mm}
		s^{2-\frac{2}{n}-\beta}
		\Big(s-\frac{1-\al}{y(t)}\Big)^{(1-\al)(\sig-m+1)+\al-2} \nn\\
	& & \hs{40mm}
	\qquad \mbox{for all $t\in (0,T)$ and $s\in (\frac{1}{y(t)},R^n) \cap (0,\ssst)$},
  \eea
  where in view of (\ref{7.13}) and (\ref{7.93}), we may utilize (\ref{7.6}) to find that
  \bas
	s^{2-\frac{2}{n}-\beta} \cdot \Big(s-\frac{1-\al}{y(t)}\Big)^{(1-\al)(\sig-m+1)+\al-2}
	&\le& c_1(\al) \cdot s^{2-\frac{2}{n}-\beta + (1-\al)(\sig-m+1)+\al-2} \\
	&=& c_1(\al) \cdot s^{-\frac{2}{n}-\beta + (1-\al)(\sig-m+1)+\al} \\
	&\le& c_1(\al) \cdot \ssst^{-\frac{2}{n}-\beta + (1-\al)(\sig-m+1)+\al} \\[2mm]
	& & \hs{10mm}
	\qquad \mbox{for all $t\in (0,T)$ and $s\in (\frac{1}{y(t)},R^n) \cap (0,\ssst)$}.
  \eas
  Furthermore estimating $e^{(\sig-m+1)\theta t} \le e^{\sig-m+1}$ for $t\in (0,\frac{1}{\theta})$, from the smallness assumption
  (\ref{7.91}) on $\ssst$ we thus infer that (\ref{7.15}) implies the inequality
  \bea{7.16}
	- n^2 s^{2-\frac{2}{n}} \KD \cdot \big( n e^{-\theta t} \hU_s \big)^{m-1} \hU_{ss}
	&\le& \frac{1}{2} \cdot \kS \cdot \big( n e^{-\theta t} \hU_s \big)^\sig \cdot \frac{a}{2} s^\beta \nn\\[2mm]
	& &
	\mbox{for all $t\in (0,T)\cap (0,\frac{1}{\theta})$ and $s\in (\frac{1}{y(t)},R^n) \cap (0,\ssst)$}.
  \eea
  We next combine (\ref{1.2}) with (\ref{1.3}) to compute
  \bea{7.17}
	\frac{\hU_t}{\kS \cdot (n e^{-\theta t} \hU_s)^\sig \cdot \frac{a}{2} s^\beta}
	&=& \frac{2 e^{\sig\theta t}}{n^\sig a \kS} \cdot s^{-\beta} \hU_s^{-\sig} \hU_t \nn\\
	&=& \frac{2 e^{\sig\theta t}}{n^\sig a \kS} \cdot s^{-\beta} \cdot
		\Big\{ \al^{1-\al} a\cdot \Big(s-\frac{1-\al}{y(t)} \Big)^{\al-1} \Big\}^{-\sig} \times \nn\\
	& & \hs{50mm}
		\al^{1-\al} (1-\al) a\cdot \Big(s-\frac{1-\al}{y(t)}\Big)^{\al-1} \cdot \frac{y'(t)}{y^2(t)} \nn\\
	&=& \frac{2 \al^{(1-\al)(1-\sig)} (1-\al) e^{\sig\theta t}}{n^\sig a^\sig \kS} \cdot s^{-\beta} \cdot
		\Big(s-\frac{1-\al}{y(t)} \Big)^{(1-\al)(\sig-1)} \cdot \frac{y'(t)}{y^2(t)} \nn\\[2mm]
	& & \hs{24mm}
	\qquad \mbox{for all $t\in (0,T)$ and $s\in (\frac{1}{y(t)},R^n) \cap (0,\ssst)$},
  \eea
  where by (\ref{7.13}), (\ref{7.94}) and (\ref{7.7}), regardless of the sign of $(1-\al)(\sig-1)$ we have
  \bas
	s^{-\beta} \cdot \Big(s-\frac{1-\al}{y(t)} \Big)^{(1-\al)(\sig-1)}
	&\le& c_2(\al) s^{-\beta+(1-\al)(\sig-1)} \\
	&\le& c_2(\al) y^{[\beta-(1-\al)(\sig-1)]_+}(t) \\
	&\le& c_2(\al) y^{1-\delss}(t)
	\qquad \mbox{for all $t\in (0,T)$ and $s\in (\frac{1}{y(t)},R^n) \cap (0,\ssst)$},
  \eas
  because $\ssst\le 1$, and because $y\ge 1$ due to (\ref{7.2}).
  On the right-hand side of (\ref{7.17}), according to (\ref{7.2}) we thus see that
  \bas
	s^{-\beta} \cdot \Big(s-\frac{1-\al}{y(t)} \Big)^{(1-\al)(\sig-1)} \cdot \frac{y'(t)}{y^2(t)}
	&\le& c_2(\al) y^{-1-\delss}(t) y'(t) \nn\\
	&\le& c_2(\al) \gamss
	\qquad \mbox{for all $t\in (0,T)$ and $s\in (\frac{1}{y(t)},R^n) \cap (0,\ssst)$},
  \eas
  so that since $e^{\sig\theta t} \le e^\sig$ for all $t\in (0,\frac{1}{\theta})$, as a consequence of (\ref{7.92}) we conclude from
  (\ref{7.17}) that
  \bas
	\hU_t
	\le \frac{1}{2} \cdot \kS \cdot \big( n e^{-\theta t} \hU_s)^\sig \cdot \frac{a}{2} s^\beta
	\qquad \mbox{for all $t\in (0,T)\cap (0,\frac{1}{\theta})$ and $s\in (\frac{1}{y(t)},R^n) \cap (0,\ssst)$}.
  \eas
  Together with (\ref{7.16}) inserted into (\ref{7.14}), this yields (\ref{7.4}).
\qed
Fortunately, whatever value of $s_0$ will eventually be chosen, in the very outer part where $s_0\le s<R^n$
a favorable subsolution property of $\uU$ can be achieved upon selecting the parameter $\theta$ in (\ref{uU})-(\ref{uW})
suitably large. This will result from the following observation.
\begin{lem}\label{lem8}
  Let $n\ge 1$ and $R>0$, assume (\ref{kl}) and (\ref{DSreg}), and let $\mus>0,\muS>0$, $\al\in (0,1)$
  and $s_0\in (0,1]$ be given.
  Then there exists $\ths=\ths^{(\mus,\muS,\al,s_0)}>0$ such that whenever $\theta>\ths$ and $T>0$ as well as $y\in C^1([0,T))$
  are such that
  \be{8.1}
	y(t)>\frac{1}{R^n}
	\quad \mbox{and} \quad
	 0\le y'(t) \le y^2(t)
	\qquad \mbox{for all } t\in (0,T),
  \ee
  the function $\hU$ from (\ref{hU}) satisfies
  \bea{8.2}
	& & \hs{-32mm}
	\hU_t - \theta \hU - n^2 s^{2-\frac{2}{n}} D\big( n e^{-\theta t} \hU_s\big) \hU_{ss}
	+ S\big( n e^{-\theta t} \hU_s\big) \cdot \frac{\muS e^{\theta t} s}{n} + \ko \hU \nn\\[2mm]
	&\le& 0
	\qquad \mbox{for all $t\in (0,T) \cap (0,\frac{1}{\theta})$ and $s\in (\frac{1}{y(t)},R^n) \cap [s_0,R^n)$.}
  \eea
\end{lem}
\proof
  Taking $a=a^{(\mus)}$ as in (\ref{a}), we use (\ref{DSreg}) to find $c_i=c_i(\mus,s_0)>0$, $i\in\{1,2\}$, such that
  \be{8.3}
	D(\xi)\le c_1
	\quad \mbox{and} \quad
	S(\xi) \le c_2
	\qquad \mbox{for all } \xi\in \Big(0,\frac{n a}{s_0}\Big),
  \ee
  and given $\al\in (0,1)$ we can then choose $\ths=\ths^{(\mus,\muS,\al,s_0)}>0$ large enough such that $\ths>\ko$ and
  \be{8.4}
	(\ths-\ko) \cdot a s_0^\al
	\ge \frac{a}{s_0}
	+ \frac{n^2 c_1 a R^{2n-2}}{\al s_0^2}
	+ \frac{c_2 \muS e R^n}{n}.
  \ee
  If $\theta>\ths$, and if $T>0$ and $y\in C^1([0,T))$ satisfy (\ref{8.1}), since
  \bas
	s-\frac{1-\al}{y(t)} \ge s-(1-\al)s=\al s
	\ge \al s_0
	\qquad \mbox{for all $t\in (0,T)$ and $s\in (\frac{1}{y(t)},R^n) \cap [s_0,R^n)$},
  \eas
  we can then use that $\al\in (0,1)$ and $s_0\le 1$ to see that in (\ref{hU}), (\ref{1.2}), (\ref{1.3}) and (\ref{1.4}),
  \bas
	\hU = \al^{-\al} a \cdot \Big(s-\frac{1-\al}{y(t)}\Big)^\al
	\ge \al^{-\al} a\cdot (\al s_0)^\al
	= a s_0^\al
  \eas
  and
  \bas
	\hU_t
	&=& \al^{1-\al} (1-\al) a \cdot \Big(s-\frac{1-\al}{y(t)}\Big)^{\al-1} \cdot \frac{y'(t)}{y^2(t)} \\
	&\le& \al^{1-\al} (1-\al) a \cdot (\al s_0)^{\al-1}
	= (1-\al) a s_0^{\al-1}
	\le \frac{a}{s_0}
  \eas
  as well as
  \bas
	\hU_s
	= \al^{1-\al} a \cdot \Big(s-\frac{1-\al}{y(t)}\Big)^{\al-1}
	\le \al^{1-\al} a \cdot (\al s_0)^{\al-1}
	= a s_0^{\al-1}
	\le \frac{a}{s_0}
  \eas
  and
  \bas
	-\hU_{ss} = \al^{1-\al} (1-\al) a \cdot \Big(s-\frac{1-\al}{y(t)}\Big)^{\al-2}
	\le \al^{1-\al} (1-\al) a \cdot (\al s_0)^{\al-2}
	= \frac{(1-\al) a s_0^{\al-2}}{\al}
	\le \frac{a}{\al s_0^2}
  \eas
  for all $t\in (0,T)$ and $s\in (\frac{1}{y(t)},R^n) \cap [s_0,R^n)$.
  Consequently, (\ref{8.3}) applies so as to ensure that
  \bas
	& & \hs{-10mm}
	\hU_t - \theta \hU - n^2 s^{2-\frac{2}{n}} D\big( n e^{-\theta t} \hU_s\big) \hU_{ss}
	+ S\big( n e^{-\theta t} \hU_s\big) \cdot \frac{\muS e^{\theta t} s}{n} + \ko \hU \\
	&\le& \frac{a}{s_0}
	- (\theta-\ko) \cdot a s_0^\al
	+ n^2 \cdot (R^n)^{2-\frac{2}{n}} \cdot c_1 \cdot \frac{a}{\al s_0^2}
	+ c_2 \cdot \frac{\muS e R^n}{n} \\[2mm]
	& & \hs{10mm}
	\qquad \mbox{for all $t\in (0,T) \cap (0,\frac{1}{\theta})$ and $s\in (\frac{1}{y(t)},R^n) \cap [s_0,R^n)$.}
  \eas
  An application of (\ref{8.4}) thus completes the proof.
\qed
Next, in an intermediate region where $\frac{1}{y(t)}<s<\sssst$ with some adequately small $\sssst$,
the action of $\qarab$ can be controlled whenever $\al$ is suitably small relative to $\beta$.
\begin{lem}\label{lem9}
  Let $n\ge 3$ and $R>0$, and assume (\ref{kl}).
  Then whenever $\mus>0, \muS>0, \beta\in (\frac{2}{n},1)$ and $\al\in (0,\beta-\frac{2}{n})$, one can find
  $\sssst=\sssst(\mus,\muS,\al,\beta)\in (0,1]$ with the property that if $T>0$ and $y\in C^1([0,T))$ are such that
  \be{9.1}
	y(t)>\max \Big\{ 1 \, , \, \frac{1}{R^n}\Big\}
	\quad \mbox{and} \quad
	0\le y'(t) \le y^{1+\frac{2}{n}}(t)
	\qquad \mbox{for all } t\in (0,T),
  \ee
  and if
  \be{9.2}
	\theta\ge \lo,
  \ee
  then the functions $\uU=\uU^{(\mus,\al,\theta,y)}$ and
  $\uW=\uW^{(\mus,\beta,\theta,y)}$ defined in (\ref{uU}) and (\ref{uW}) satisfy
  \be{9.3}
	\qarab [\uU,\uW](s,t) \le 0
	\qquad \mbox{for all $t\in (0,T)$ and $s\in (\frac{1}{y(t)},R^n) \cap (0,\sssst)$.}
  \ee
\end{lem}
\proof
  We use that $\beta-\frac{2}{n}-\al>0$ in choosing
  $\sssst=\sssst(\mus,\muS,\al,\beta)\in (0,1]$ such that
  \be{9.4}
	\frac{n^2}{\lt \beta} \sssst^{\beta-\frac{2}{n}-\al} \le \frac{1}{2}
  \ee
  as well as
  \be{9.5}
	\frac{1}{\lt} \sssst^{\beta-\frac{2}{n}-\al} \le \frac{1}{2}.
  \ee
  Then assuming that $T>0$ and $y\in C^1([0,T))$ are such that (\ref{9.1}) holds, for $\theta\ge \lo$ we obtain from
  (\ref{Q}), (\ref{hU}), (\ref{2.2}), (\ref{2.4}) and (\ref{9.1}) that
  \bea{9.6}
	e^{\theta t} \qarab [\uU,\uW](s,t)
	&=& \hW_t - \theta \hW - n^2 s^{2-\frac{2}{n}} \hW_{ss} + \lo \hW - \lt \hU \nn\\
	&=& \beta^{1-\beta} (1-\beta) a \cdot \Big( s-\frac{1-\beta}{y(t)}\Big)^{\beta-1} \cdot \frac{y'(t)}{y^2(t)} \nn\\
	& & + n^2 s^{2-\frac{2}{n}} \cdot \beta^{1-\beta} (1-\beta) a \cdot \Big( s-\frac{1-\beta}{y(t)}\Big)^{\beta-2} \nn\\
	& & - \lt \cdot \al^{-\al} a \cdot \Big( s- \frac{1-\al}{y(t)} \Big)^\al \nn\\
	& & - (\theta-\lo)\hW \nn\\
	&\le& \beta^{1-\beta} (1-\beta) a \cdot \Big( s-\frac{1-\beta}{y(t)}\Big)^{\beta-1} \cdot y^{-1+\frac{2}{n}}(t) \nn\\
	& & + n^2 \beta^{1-\beta} a s^{2-\frac{2}{n}} \cdot \Big( s-\frac{1-\beta}{y(t)}\Big)^{\beta-2} \nn\\
	& & - \lt \al^{-\al} a \cdot \Big( s- \frac{1-\al}{y(t)} \Big)^\al
	\qquad \mbox{for all $t\in (0,T)$ and $s\in (\frac{1}{y(t)},R^n)$},
  \eea
  because $0\le 1-\beta \le 1$ and $\hW\ge 0$.
  Here since
  \bas
	s-\frac{1-\al}{y(t)} \ge \al s
	\quad \mbox{and} \quad
	s-\frac{1-\beta}{y(t)} \ge \beta s
	\qquad \mbox{for all $t\in (0,T)$ and $s\in (\frac{1}{y(t)},R^n)$},
  \eas
  from (\ref{9.4}) it follows that
  \bas
	\frac{n^2 \beta^{1-\beta} a s^{2-\frac{2}{n}} (s-\frac{1-\beta}{y(t)})^{\beta-2}}
		{\lt \al^{-\al} a (s-\frac{1-\al}{y(t)})^\al}
	&\le& \frac{n^2 \beta^{1-\beta} a s^{2-\frac{2}{n}} \cdot (\beta s)^{\beta-2}}
		{\lt \al^{-\al} a \cdot (\al s)^\al} \\
	&=& \frac{n^2}{\lt \beta} \cdot s^{\beta-\frac{2}{n}-\al} \\
	&\le& \frac{1}{2}
	\qquad \mbox{for all $t\in (0,T)$ and $s\in (\frac{1}{y(t)},R^n) \cap (0,\sssst)$},
  \eas
  while (\ref{9.5}) similarly implies that
  \bas
	\frac{\beta^{1-\beta} a (s-\frac{1-\beta}{y(t)})^{\beta-1} y^{-1+\frac{2}{n}}(t)}
		{\lt \al^{-\al} a (s-\frac{1-\al}{y(t)})^\al}
	&\le& \frac{s^{\beta-1} y^{-1+\frac{2}{n}}(t)}{\lt s^\al} \\
	&=& \frac{1}{\lt} \cdot s^{\beta-1 -\al} y^{-1+\frac{2}{n}}(t) \\
	&\le& \frac{1}{\lt} \cdot s^{\beta-1 -\al} \cdot s^{1-\frac{2}{n}} \\
	&=& \frac{1}{\lt} \cdot s^{\beta-\frac{2}{n}-\al} \\
	&\le& \frac{1}{2}
	\qquad \mbox{for all $t\in (0,T)$ and $s\in (\frac{1}{y(t)},R^n) \cap (0,\sssst)$}.
  \eas
  From (\ref{9.6}) we therefore obtain (\ref{9.3}).
\qed
Again for arbitrary choices of $s_0$, a corresponding feature near the outer boundary can be asserted by taking
$\theta$ large enough:
\begin{lem}\label{lem10}
  Let $n\ge 3$ and $R>0$, and assume (\ref{kl}).
  Then for each $\beta\in (\frac{2}{n},1)$ and any $s_0\in (0,1]$ there exists $\thss=\thss^{(\beta,s_0)}>0$ such that if
  $\mus>0, \al\in (0,1)$,
  $\theta>\thss, T>0$ and $y\in C^1([0,T))$ are such that
  \be{10.1}
	y(t)>\max\Big\{ 1\, , \, \frac{1}{R^n}\Big\}
	\quad \mbox{and} \quad
	 0\le y'(t) \le \lt y^{1-\al+\beta}(t)
	\qquad \mbox{for all } t\in (0,T),
  \ee
  it follows that with $\uU=\uU^{(\mus,\al,\theta,y)}$ and
  $\uW=\uW^{(\mus,\beta,\theta,y)}$ taken from (\ref{uU}) and (\ref{uW}),
  \be{10.2}
	\qarab [\uU,\uW](s,t) \le 0
	\qquad \mbox{for all $t\in (0,T)$ and $s\in (\frac{1}{y(t)},R^n) \cap [s_0,R^n)$.}
  \ee
\end{lem}
\proof
  We let $\thss=\thss^{(\beta,s_0)}>0$ be such that
  \be{10.3}
	\thss\ge 2\lo
	\qquad \mbox{and} \qquad
	\thss \ge \frac{2n^2}{\beta s_0^\frac{2}{n}},
  \ee
  and then observe that if $\mus>0$, $\al\in (0,1)$, $\theta>\thss$, $T>0$ and $y\in C^1([0,T))$ are such that (\ref{10.1})
  holds, in (\ref{Q}) we can estimate
  \bas
	e^{\theta t} \qarab [\uU,\uW](s,t)
	&=& \hW_t - n^2 s^{2-\frac{2}{n}} \hW_{ss} - \lt \hU - (\theta-\lo)\hW \\
	&\le& \hW_t - n^2 s^{2-\frac{2}{n}} \hW_{ss} - \lt \hU - \frac{\theta}{2} \hW
	\qquad \mbox{for all $t\in (0,T)$ and $s\in (\frac{1}{y(t)},R^n)$.}
  \eas
  Again since $s-\frac{1-\beta}{y(t)} \ge \beta s$ and $s-\frac{1-\al}{y(t)}\ge \al s$
  for all $t\in (0,T)$ and $s\in (\frac{1}{y(t)},R^n)$, in line with (\ref{10.1}), Lemma \ref{lem2} and Lemma \ref{lem1} this
  implies that for all $t\in (0,T)$ and $s\in (\frac{1}{y(t)},R^n)$,
  \bea{10.4}
	e^{\theta t} \qarab [\uU,\uW](s,t)
	&\le& \beta^{1-\beta} (1-\beta) a\cdot\Big( s-\frac{1-\beta}{y(t)}\Big)^{\beta-1} \cdot \lt y^{\beta-\al-1}(t) \nn\\
	& & + n^2 s^{2-\frac{2}{n}} \cdot \beta^{1-\beta} (1-\beta) a\cdot \Big(s-\frac{1-\beta}{y(t)}\Big)^{\beta-2} \nn\\
	& & - \lt \al^{-\al} a \cdot \Big( s-\frac{1-\al}{y(t)}\Big)^\al \nn\\
	& & - \frac{\theta}{2} \cdot \beta^{-\beta} a\cdot\Big( s-\frac{1-\beta}{y(t)}\Big)^\beta \nn\\
	&\le& \lt a s^{\beta-1} y^{\beta-\al-1}(t)
	+ \frac{n^2 a}{\beta} \cdot s^{\beta-\frac{2}{n}}
	- \lt a s^\al
	- \frac{\theta a}{2} \cdot s^\beta,
  \eea
  because $0<\beta<1$.
  But here the inequality $\beta-\al-1<0$ warrants that
  \bas
	\lt a s^{\beta-1} y^{\beta-\al-1}(t)
	\le \lt a s^{\beta-1} \cdot \Big(\frac{1}{s}\Big)^{\beta-\al-1}
	= \lt a s^\al
	\qquad \mbox{for all $t\in (0,T)$ and $s\in (\frac{1}{y(t)},R^n)$,}
  \eas
  while utilizing the second restriction in (\ref{10.3}) shows that in the considered region we have
  \bas
	\frac{\frac{n^2 a}{\beta} \cdot s^{\beta-\frac{2}{n}}}{\frac{\theta a}{2} \cdot s^\beta}
	&=& \frac{2n^2}{\beta\theta} \cdot s^{-\frac{2}{n}}
	\le \frac{2n^2}{\beta\theta} \cdot s_0^{-\frac{2}{n}} \nn\\
	&\le& 1
	\qquad \mbox{for all $t\in (0,T)$ and $s\in (\frac{1}{y(t)},R^n) \cap [s_0,R^n)$.}
  \eas
  Therefore, (\ref{10.2}) results from (\ref{10.4}).
\qed
\subsection{Conclusion. Proof of Theorem \ref{theo13}}\label{sect_BU_proof}
We can now make sure that under the assumptions on $\sig$ from Theorem \ref{theo13} we can find
some small $\al\in (0,1)$, some $\beta\in (\frac{2}{n},1)$ close to $\frac{2}{n}$, some large $\theta>0$ and some
positive function $y$, blowing up at some suitably small $T>0$, such that
all the hypotheses in Lemmata \ref{lem5}-\ref{lem10} are simultaneously fulfilled,
and that hence the above construction indeed yields a subsolution pair $(\uU,\uW)$ for which $\uU$ develops
a shock-type singularity at $(s,t)=(0,T)$.
\begin{lem}\label{lem12}
  Let $n\ge 3$ and $R>0$, assume (\ref{kl}), and suppose that (\ref{DSreg}), (\ref{Du}) and (\ref{Sl}) hold with some
  $\xi_0>0, \KD>0, \kS>0, m\in\R$ and $\sig\in\R$ such that
  \be{12.1}
	\sig>m-1+\frac{4}{n}
	\qquad \mbox{and} \qquad
	\sig>\frac{2}{n}.
  \ee
  Then there exist $\al\in (0,1)$ and $\beta\in (\frac{2}{n},1)$ such that whenever $\mus>0, \muS>0$ and $\Ts>0$,
  one can find $\theta>0, T\in (0,\Ts)$ and a positive function $y\in C^1([0,T))$ such that
  \be{12.2}
	y(t) \to +\infty
	\qquad \mbox{as } t\nearrow T,
  \ee
  and that with the functions $\uU=\uU^{(\mus,\al,\theta,y)}$ and
  $\uW=\uW^{(\mus,\beta,\theta,y)}$ defined through (\ref{uU}) and (\ref{uW}) satisfy
  \be{12.3}
	\parab^{(\muS)} [\uU,\uW](s,t) \le 0
	\quad \mbox{and} \quad
	\qarab [\uU,\uW](s,t) \le 0
	\qquad \mbox{for all $t\in (0,T)$ and $s\in (0,R^n)\sm \{ \frac{1}{y(t)} \}$.}
  \ee
\end{lem}
\proof
  With $\als,\bets,\dels,\alss,\betss$ and $\delss$ as provided by Lemma \ref{lem5} and Lemma \ref{lem7}, we pick any
  $\al\in (0,1)$ and $\beta\in (\frac{2}{n},1)$ such that
  \be{12.4}
	\al<\als,
	\quad
	\al<\alss,
	\quad
	\beta<\bets,
	\quad
	\beta<\betss
	\quad \mbox{and} \quad
	\al<\beta-\frac{2}{n},
  \ee
  and let
  \be{12.5}
	\del:=\min \Big\{ \dels \, , \, \delss \, , \, \beta-\al \, , \, \frac{2}{n} \Big\}.
  \ee
  Given $\mus>0, \muS>0$ and $\Ts>0$, we then let $\ssst=\ssst^{(\mus,\muS)}$ and $\sssst=\sssst^{(\mus,\muS,\al,\beta)}$
  be as in Lemma \ref{lem7} and Lemma \ref{lem9}, and apply Lemma \ref{lem8} and Lemma \ref{lem10} to
  \be{12.6}
	s_0:=\min \big\{ \ssst \, , \, \sssst \big\}
  \ee
  to find $\ths=\ths^{(\mus,\muS,\al,s_0)}$ and $\thss=\thss^{(\beta,s_0)}$ with the properties listed there. Letting
  \be{12.7}
	\theta:=\max \big\{ \ths \, , \, \thss \, , \, \ko \, , \, \lo \big\}
  \ee
  and taking $a=a^{(\mus)}$, $\gams=\gams^{(\mus,\muS)}$, $\ys=\ys^{(\mus,\muS)}$ and $\gamss=\gamss^{(\mus,\muS,\al,\beta)}$
  from (\ref{a}), lemma \ref{lem5} and Lemma \ref{lem7}, we then set
  \be{12.8}
	\gamma := \min \Big\{ \gams \, , \, \gamss \, , \, \lt \, ,\, 1\Big\}
  \ee
  and choose $y_0>0$ large enough such that
  \be{12.9}
	y_0>\max \Big\{ 1 \, \frac{1}{R^n}\Big\}
	\qquad \mbox{and} \qquad
	y_0 \ge \ys,
  \ee
  and that
  \be{12.10}
	T:=\frac{1}{\gamma\del y_0^\del}
	\quad \mbox{satisfies} \quad
	T < \min \Big\{ \frac{1}{\theta} \, , \, \Ts \Big\}.
  \ee
  By explicit computation, this choice of $T$ can be seen to ensure that the problem
  \be{12.11}
	\lbal
	y'(t) = \gamma y^{1+\del}(t),
	\qquad t\in (0,T), \\[1mm]
	y(0)=y_0,
	\ear
  \ee
  possesses a solution $y\in C^1([0,T))$ which satisfies (\ref{12.2}) as well as $y\ge y_0$, so that (\ref{12.9}),
  (\ref{12.8}) and (\ref{12.5}) guarantee that all the requirements on $y$ made in Lemmata \ref{lem5}-\ref{lem10} are met.
  As furthermore (\ref{12.4}) and (\ref{12.7}) warrant that also the corresponding assumptions on $\al$, $\beta$ and $\theta$
  are fulfilled, we may draw on said statements and on (\ref{12.6}) to readily infer that for
  $(\uU,\uW)=(\uU^{(\mus,\al,\theta,y)},\uW^{(\mus,\beta,\theta,y)})$
  as accordingly defined through (\ref{uU}) and (\ref{uW}) we have
  \bas
	\parab^{(\muS)} [\uU,\uW](s,t)\le 0
	\quad \mbox{and} \quad
	\qarab [\uU,\uW](s,t)\le 0
	\qquad \mbox{for all $t\in (0,T) \cap (0,\frac{1}{\theta})$ and $s\in (0,R^n)\sm \{ \frac{1}{y(t)}\}$.}
  \eas
  Since from (\ref{12.10}) we moreover know that $(0,T)\cap (0,\frac{1}{\theta})=(0,T)$ and that $T<\Ts$, the proof thereby
  becomes complete.
\qed
Thanks to the comparison principle from Lemma \ref{lem_CP}, in the above constellation any solution emanating from initial data
for which $U(\cdot,0)\ge \uU(\cdot,0)$ and $W(\cdot,0)\ge \uW(\cdot,0)$ must cease to exist within finite time:\abs
\proofc of Theorem \ref{theo13}. \quad
  Given $\TS>0, \Ms>0$ and $\MS>\Ms$, with $T_0=T_0^{(\Ms,\MS)}$ taken from Lemma \ref{lem11} we let
  \be{13.4}
	\Ts:= \min \big\{ T_0 \, , \, \TS \big\}
  \ee
  and apply Lemma \ref{lem12} to $\mus$ and $\muS$ as in (\ref{11.5}) to find
  $\al\in (0,1)$ and $\beta\in (\frac{2}{n},1)$, $\theta>0, T\in (0,\Ts)$ and a positive
  $y\in C^1([0,T))$ such that (\ref{12.2}) and (\ref{12.3}) hold for the pair
  $(\uU,\uW)=(\uU^{(\mus,\al,\theta,y)},\uW^{(\mus,\beta,\theta,y)})$
  correspondingly introduced in (\ref{uU}) and (\ref{uW}).
  Setting
  \be{13.5}
	M^{(u)}(r):=n|B_1(0)| \uU(r^n,0)
	\quad \mbox{and} \quad
	M^{(w)}(r):=n|B_1(0)| \uW(r^n,0)
	\qquad \mbox{for } r\in [0,R]
  \ee
  then defines continuous functions $M^{(u)}$ and $M^{(w)}$ on $[0,R]$ which due to the inclusions
  $\uU(\cdot,0)\in C^1([0,R^n])$ and $\uW(\cdot,0)\in C^1([0,R^n])$ as well as the identities $\uU(0,0)=0$ and $\uW(0,0)=0$
  satisfy (\ref{13.01}), and for which thanks to Lemma \ref{lem1} and Lemma \ref{lem2} we furthermore have
  $M^{(u)}(R) \le n|B_1(0)| \cdot \frac{\mus R^n}{n}=\mus |\Om|$
  and, similarly, $M^{(w)}(R)\le \mus |\Om|$, whence also (\ref{13.1}) holds due to (\ref{11.5}).\abs
  Now if $u_0$ and $w_0$ satisfy (\ref{ir}), (\ref{13.11}) and (\ref{13.2}), then assuming that
  for the solution $(u,v,w)$ of (\ref{0})
  from Lemma \ref{lem_loc} we had $\tm\ge \TS$ we could draw on (\ref{13.11}) and the inequality $\Ts\le T_0$, as implied
  by (\ref{13.4}), to employ Lemma \ref{lem11} to thus conclude that the functions $U$ and $W$ from (\ref{U}) and (\ref{W})
  would satisfy
  \be{13.6}
	\parab^{(\muS)} [U,W](s,t)   \ge  0
	\quad \mbox{and} \quad
	\qarab [U,W](s,t) = 0
	\qquad \mbox{for all $s\in 0,R^n)$ and } t\in (0,\Ts)
  \ee
  as well as
  \be{13.7}
	U(0,t)=W(0,t)=0,
	\quad U(R^n,t) \ge \frac{\mus R^n}{n}
	\quad \mbox{and} \quad
	W(R^n,t) \ge \frac{\mus R^n}{n}
	\qquad \mbox{for all } t\in (0,\Ts).
  \ee
  On the other hand, (\ref{uU}) and (\ref{uW}) together with (\ref{1.1}) and (\ref{2.1}) assert that
  for all $t\in (0,T)$,
  \bas
	\uU(0,t)=\uW(0,t)=0,
	\quad
	\uU(R^n,t) = e^{-\theta t} \hU(R^n,t) \le \frac{\mus R^n}{n}
	\quad \mbox{and} \quad
	\uW(R^n,t) = e^{-\theta t} \hW(R^n,t) \le \frac{\mus R^n}{n}.
  \eas
  Since furthermore a combination of (\ref{13.2}) with (\ref{13.5}) shows that
  \bas
	U(s,0)=\frac{1}{n|B_1(0)|} \int_{B_{s^{1/n}}(0)} u_0 dx \ge \frac{1}{n|B_1(0)|} M^{(u)}(s^\frac{1}{n}) = \uU(s,0)
	\qquad \mbox{for all } s\in (0,R^n),
  \eas
  and that likewise also
  \bas
	W(s,0)\ge \uW(s,0)
	\qquad \mbox{for all } s\in (0,R^n),
  \eas
  relying on (\ref{13.6}) and (\ref{12.3}) we may invoke the comparison principle from Lemma \ref{lem_CP} to infer that
  \bas
	U(s,t) \ge \uU(s,t)
	\qquad \mbox{for all $s\in (0,R^n)$ and } t\in (0,T).
  \eas
  In view of (\ref{1.3}), this would particularly imply that with $a=a^{(\mus)}$ as in (\ref{a}),
  \bas
	\frac{1}{n} \cdot u(0,t)
	= U_s(0,t)
	\ge \uU_s(0,t)
	=  e^{-\theta t}\cdot ay^{1-\al}(t) \ge a e^{-\theta T}\cdot y^{1-\al}(t)
	\qquad \mbox{for all } t\in (0,T),
  \eas
  which according to (\ref{12.2}) is incompatible with our hypothesis that $\tm\ge \Ts$, because $T<\Ts$.
\qed
\mysection{Global existence for subcritical $\sig$} \label{sect_4}
This section will reveal optimality of the condition (\ref{sig}) by asserting global solvablity, and partially also boundedness,
under assumptions on $\sig$ which are essentially complementary to those from Theorem \ref{theo13}.
Our analysis in this direction pursues an approach of basically variational nature, operating at the level of the original
variables and not relying on any symmetry requirements.
The starting point will be formed by the following outcome of a standard testing procedure applied to the first
equation in (\ref{0}).
\begin{lem}\label{lem20}
  Suppose that $n\ge 3$, that $\Om\subset\R^n$ is a bounded domain with smooth boundary,
  that (\ref{kl}), (\ref{DSreg}), (\ref{Dl}), (\ref{Su}) and (\ref{init}) hold with some $\xi_0>0, \kD>0, \KS>0, m\in\R$ and
  $\sig\in\R$.
  Then for each $p>\max\{1,1-\sig\}$, one can find $C(p)>0$ such that
  \bea{20.1}
	& & \hs{-20mm}
	\frac{d}{dt} \io (u+1)^p dx
	+ \frac{1}{C(p)} \io \big|\na (u+1)^\frac{p+m-1}{2}\big|^2 dx
	+ p\ko \io (u+1)^{p-1} u dx \nn\\
	&\le& C(p) \io (u+1)^{p+\sig-1} w dx
	+ p\kt \io (u+1)^{p-1} w dx
	\qquad \mbox{for all } t\in (0,\tm).
  \eea
\end{lem}
\proof
  According to an integration by parts in the first equation from (\ref{0}),
  \bea{20.2}
	\frac{1}{p} \frac{d}{dt} \io (u+1)^p dx
	&=& - (p-1) \io (u+1)^{p-2} D(u) |\na u|^2 dx
	+ (p-1) \io (u+1)^{p-2} S(u) \na u\cdot \na v dx \nn\\
	& & - \ko \io (u+1)^{p-1} u dx
	+ \kt \io (u+1)^{p-1} w dx
	\qquad \mbox{for all } t\in (0,\tm),
  \eea
  where a second integration by parts using the second equation in (\ref{0}) shows that
  for all $t\in (0,\tm)$.
  \bea{20.3}
	(p-1) \io (u+1)^{p-2} S(u) \na u \cdot\na v dx
	&=& \io \na\psi_p(u)\cdot\na v dx
	= - \io \psi_p(u)\Del v dx \nn\\
	&=& \io \psi_p(u) w dx - \muw(t) \io \psi_p(u) dx
  \eea
  with
  \bas
	\psi_p(\xi):=(p-1)\int_0^\xi (\tau +1)^{p-2} S(\tau) d\tau,
	\qquad \xi\ge 0.
  \eas
  Since (\ref{DSreg}), (\ref{Dl}) and (\ref{Su}) along with our assumptions $p>1$ and $p>1-\sig$ readily imply the existence
  of $c_1(p)>0$ and $c_2(p)>0$ such that
  \bas
	(p-1) (\xi+1)^{p-2} D(\xi) \ge c_1(p) (\xi+1)^{p+m-3}
	\quad \mbox{and} \quad
	\psi_p(\xi) \le c_2(p)(\xi+1)^{p+\sig-1}
	\qquad \mbox{for all } \xi\ge 0,
  \eas
  due to the nonnegativity of $\ko$ and $\muw$ we obtain (\ref{20.1}) from (\ref{20.2}) and (\ref{20.3}) if we let
  $C(p):=\max \Big\{ \frac{(p+m-1)^2}{4p c_1(p)} \, , \, p c_2(p) \Big\}$, for instance.
\qed
The two summands in (\ref{20.1}) which involve $w$ will be estimated by means of the following basically well-known
consequence of $L^p$-$L^q$ regulatization properties of the Neumann heat semigroup.
In dependence of the particular situations addressed in Proposition \ref{prop23} and Proposition \ref{prop25},
the applications of this in Lemma \ref{lem22} and in Lemma \ref{lem24} will rely on different choices of the
parameters $k$ and $l$ appearing herein.
\begin{lem}\label{lem21}
  Let $n\ge 3$ and $\Om\subset\R^n$ be a bounded domain with smooth boundary,
  and assume (\ref{kl}), (\ref{DSreg}) and (\ref{init}).
  Then for all $k\ge 1$ and $l\in [1,\frac{nk}{(n-2k)_+})$ there exists $C(k,l)>0$ such that
  \be{21.1}
	\|w(\cdot,t)\|_{L^l(\Om)} \le C(k,l) \cdot \Big\{ 1 + \sup_{t'\in (0,t)} \|u(\cdot,t')\|_{L^k(\Om)} \Big\}
	\qquad \mbox{for all } t\in (0,\tm).
  \ee
\end{lem}
\proof
  This can be seen by applying known smoothing properties of the Neumann heat semigroup to the third equation in (\ref{0}).
\qed
\subsection{Criticality of the line $\sig=m-1+\frac{4}{n}$. Proof of Proposition \ref{prop23}}
We next focus on the parameter setting of Proposition \ref{prop23}, in which the influence of chemotaxis is assumed to be suitably small
relative to diffusion, and may then rely on the dissipative effect generated by the second summand on the left of (\ref{20.1})
to control ill-signed contributions containing $u$ by means of interpolation using the zero-order information
from (\ref{mass_rough}) and (\ref{mass_kl2}), respectively. This will be prepared by the following statement
derived from the Gagliardo-Nirenberg inequality in a straightforward manner.
\begin{lem}\label{lem211}
  Let $n\ge 3, m\in\R, p>2-m-\frac{2}{n}$ and $r>1$ be such that
  \be{211.1}
	r < \frac{n(p+m-1)}{n-2}.
  \ee
  Then there exists $C>0$ such that
  \be{211.2}
	\bigg\{ \io \vp^r dx \bigg\}^\frac{n(p+m-1) +2-n}{n(r-1)}
	\le C\cdot \bigg\{ \io \vp dx \bigg\}^\frac{n(p+m-1)-(n-2)r}{n(r-1)} \cdot
	\io \big|\na\vp^\frac{p+m-1}{2}\big|^2 dx
	+ C\cdot \bigg\{ \io \vp dx \bigg\}^\frac{[n(p+m-1)+2-n]\cdot r}{n(r-1)}
  \ee
  for all $\vp\in C^1(\bom)$ fulfilling $\vp>0$ in $\bom$.
\end{lem}
\proof
  Since $p>2-m-\frac{2}{n}>1-m$, and since
  $\frac{2}{p+m-1}<\frac{2r}{p+m-1} \le \frac{2n}{n-2}$ according to (\ref{211.1}),
  the Gagliardo-Nirenberg inequality applies so as to yield $c_1>0$ such that for all positive $\vp\in C^1(\bom)$,
  \bea{211.3}
	& & \hs{-20mm}
	\bigg\{ \io \vp^r dx \bigg\}^\frac{n(p+m-1)+2-n}{n(r-1)} \nn\\
	&=& \big\| \vp^\frac{p+m-1}{2} \big\|_{L^\frac{2r}{p+m-1}(\Om)}^{\frac{2r}{p+m-1}\cdot\frac{n(p+m-1)+2-n}{n(r-1)}} \nn\\
	&\le& c_1\big\| \na \vp^\frac{p+m-1}{2} \big\|_{L^2(\Om)}^{\frac{2r}{p+m-1}\cdot\frac{n(p+m-1)+2-n}{n(r-1)}\cdot a}
	\cdot \big\| \vp^\frac{p+m-1}{2} \big\|_{L^\frac{2}{p+m-1}(\Om)}^
		{\frac{2r}{p+m-1}\cdot \frac{n(p+m-1)+2-n}{n(r-1)} \cdot (1-a)} \nn\\
	& & + c_1 \big\| \vp^\frac{p+m-1}{2} \big\|_{L^\frac{2}{p+m-1}(\Om)}^{\frac{2r}{p+m-1}\cdot\frac{n(p+m-1)+2-n}{n(r-1)}},
  \eea
  where
  \bas
	a:=\frac{n(p+m-1)}{n(p+m-1)+2-n} \cdot \frac{r-1}{r} \ \in (0,1).
  \eas
  As thus
  \bas
	\frac{2r}{p+m-1}\cdot\frac{n(p+m-1)+2-n}{n(r-1)}\cdot a =2
  \eas
  and
  \bas
	\frac{2r}{p+m-1}\cdot\frac{n(p+m-1)+2-n}{n(r-1)}\cdot (1-a)
	= \frac{2}{p+m-1} \cdot \frac{n(p+m-1)-(n-2)r}{n(r-1)},
  \eas
  noting that for any such $\vp$ we have
  $\|\vp^\frac{p+m-1}{2}\|_{L^\frac{2}{p+m-1}(\Om)}^\frac{2}{p+m-1}= \io \vp dx$, from (\ref{211.3}) we readily obtain (\ref{211.1}).
\qed
Now under the assumption (\ref{sig2}), for arbitrarily large finite $p$ we may apply Lemma \ref{lem211} to
$r:=p+m-1+\frac{2}{n}$, and combine Lemma \ref{lem21} with an independent energy-type argument addressing $w$,
to turn Lemma \ref{lem20} into an $L^p$ bound for $u$ in the following sense.
\begin{lem}\label{lem22}
  Let $n\ge 3$ and $\Om\subset\R^n$ be a bounded domain with smooth boundary,
  and suppose that (\ref{kl}), (\ref{DSreg}), (\ref{Dl}), (\ref{Su}) is valid with some $\xi_0>0, \kD>0, \KS>0, m\in\R$
  and $\sig\in\R$ fulfilling (\ref{sig2}).
  Then whenever (\ref{init}) holds, for all $p>1$ and any $T>0$ there exists $C(p,T)>0$ such that
  \be{22.1}
	\|u(\cdot,t)\|_{L^p(\Om)} \le C(p,T)
	\qquad \mbox{for all } t\in (0,\tm) \cap (0,T),
  \ee
  and that
  \be{22.01}
	\sup_{T>0} C(p,T)<\infty
	\qquad \mbox{if (\ref{kl2}) is satisfied.}
  \ee
\end{lem}
\proof
  In view of the H\"older inequality and Young's inequality, it is sufficient to consider the case when $p>1$ is such that
  \be{22.2}
	p+m-1>\frac{n-2}{n},
	\qquad
	p+\sig-1>0,
	\qquad
	\frac{n(p+m-1)}{n-2} \ge 2
	\qquad \mbox{and} \qquad
	p>\frac{n(1-m)}{2},
  \ee
  and that
  \be{22.3}
	\frac{p+m-1+\frac{2}{n}}{m+\frac{2}{n}-\sig} < \frac{n(p+m-1)}{n-2},
  \ee
  where the latter indeed holds for all suitably large $p$ due to the fact that
  \bas
	\frac{1}{m+\frac{2}{n}-\sig}
	< \frac{1}{m+\frac{2}{n}- (m-1+\frac{4}{n})}
	= \frac{1}{1-\frac{2}{n}}
	= \frac{n}{n-2}
  \eas
  due to (\ref{sig2}), ensuring that for sufficiently large $p$ we have
  \bas
  \Big(\frac{n}{n-2}-\frac{1}{m+\frac{2}{n}-\sig}\Big) \cdot p
  >\frac{m-1+\frac{2}{n}}{m+\frac{2}{n}-\sig}
   -\frac{n(m-1)}{n-2},
  \eas which is equivalent to (\ref{22.3}).\\
  We then employ Lemma \ref{lem20} to fix $c_1=c_1(p)>0$ and $c_2=c_2(p)>0$ such that
  \bea{22.4}
	& & \hs{-20mm}
	\frac{d}{dt} \io (u+1)^p dx
	+ c_1 \io \big|\na (u+1)^\frac{p+m-1}{2}\big|^2 dx
	+ p\ko \io (u+1)^{p-1} u dx \nn\\
	&\le& c_2 \io (u+1)^{p+\sig-1} w dx
	+ p\kt \io (u+1)^{p-1} w dx
	\qquad \mbox{for all } t\in (0,\tm),
  \eea
  and for clarity in presentation first concentrating on the case when (\ref{kl2}) holds,
  we recall from (\ref{mass_kl2}) that then with some $c_3>0$,
  \be{22.5}
	\io (u +1) dx \le c_3
	\qquad \mbox{for all } t\in (0,\tm).
  \ee
  As regardless of whether or not $\ko$ is positive we know from (\ref{kl2}) that $\kt \le c_4 \ko$ with $c_4:=\frac{\lo}{\lt}>0$,
  we may twice use Young's inequality to see that with some $c_5=c_5(p)>0$ we have
  \bas
	p\kt \io (u+1)^{p-1} w dx
	&\le& p c_4 \ko \io (u+1)^{p-1} w dx \nn\\
	&\le& \ko \io (u+1)^p dx
	+ c_5 \io w^p dx
	\qquad \mbox{for all } t\in (0,\tm),
  \eas
  and that hence
  \bas
	\io (u+1)^p dx
	&\le& \io (u+1)^{p-1} u dx
	+ (u+1)^{p-1} dx \nn\\
	&\le& \io (u+1)^{p-1} u dx
	+ \frac{p-1}{p} \io (u+1)^p dx
	+ \frac{|\Om|}{p}
	\qquad \mbox{for all } t\in (0,\tm),
  \eas
  so that, in fact,
  \be{22.6}
	p\kt \io (u+1)^{p-1} w dx
	\le p\ko \io (u+1)^{p-1} u dx
	+ \ko |\Om|
	+ c_5 \io w^p dx
	\qquad \mbox{for all } t\in (0,\tm).
  \ee
  We next use that the first restriction in (\ref{22.2}) ensures applicability of Lemma \ref{lem211} to $r:=p+m-1+\frac{2}{n}$,
  because this number accordingly satisfies $r>1$ and $(n-2)r-n(p+m-1) = - 2(p+m-1) + \frac{2(n-2)}{n} <0$.
  In line with (\ref{22.5}) and noting that $\frac{n(p+m-1)+2-n}{n(r-1)}=1$, we thereby obtain $c_6=c_6(p)>0$ and $c_7=c_7(p)>0$ such that
  \be{22.7}
	\frac{c_1}{2} \io \big| \na (u+1)^\frac{p+m-1}{2}\big|^2 dx
	\ge c_6 \io (u+1)^{p+m-1+\frac{2}{n}} dx - c_7
	\qquad \mbox{for all } t\in (0,\tm),
  \ee
  and observing that (\ref{sig2}) particularly warrants that $\sig<m+\frac{2}{n}$ due to the fact that $n\ge 3$, and that
  $p+\sig-1>0$ by (\ref{22.2}), we may again rely on Young's inequality to find $c_8=c_8(p)>0$ fulfilling
  \be{22.8}
	c_2 \io (u+1)^{p+\sig-1} w dx
	\le c_6 \io (u+1)^{p+m-1+\frac{2}{n}} dx
	+ c_8 \io w^\frac{p+m-1+\frac{2}{n}}{m+\frac{2}{n}-\sig} dx
	\qquad \mbox{for all } t\in (0,\tm).
  \ee
  Now in order to appropriately control the rightmost summands in (\ref{22.6}) and (\ref{22.8}), we let
  \be{22.87}
	q:=\frac{n(p+m-1)}{n-2}
  \ee
  and note that $q\ge 2$ by (\ref{22.2}).
  Therefore, the third equation in (\ref{0}) implies that
  \be{22.9}
	\frac{1}{q} \frac{d}{dt} \io w^q dx
	+ \frac{4(q-1)}{q^2} \io |\na w^\frac{q}{2}|^2 dx
	+ \lo \io w^q dx
	= \lt \io u w^{q-1} dx
	\qquad \mbox{for all } t\in (0,\tm),
  \ee
  where a combination of the H\"older inequality with a Sobolev embedding inequality and Young's inequality shows that with
  some $c_9=c_9(p)>0$ and $c_{10}=c_{10}(p)>0$,
  \bea{22.10}
	\lt \io u w^{q-1} dx
	&\le& \lt \cdot \bigg\{ \io u^\frac{nq}{n+2q-2} dx \bigg\}^\frac{n+2q-2}{nq} \cdot
		\bigg\{ \io w^\frac{nq}{n-2} dx \bigg\}^\frac{(n-2)(q-1)}{nq} \nn\\
	&=& \lt \cdot \bigg\{ \io u^\frac{nq}{n+2q-2} dx \bigg\}^\frac{n+2q-2}{nq} \cdot
		\|w^\frac{q}{2}\|_{L^\frac{2n}{n-2}(\Om)}^\frac{2(q-1)}{q} \nn\\
	&\le& c_9 \cdot \bigg\{ \io u^\frac{nq}{n+2q-2} dx \bigg\}^\frac{n+2q-2}{nq} \cdot
		\bigg\{ \|\na w^\frac{q}{2}\|_{L^2(\Om)}^2 +1 \bigg\}^\frac{q-1}{q} \nn\\
	&\le& \frac{4(q-1)}{q^2} \cdot \bigg\{ \io |\na w^\frac{q}{2}|^2 dx + 1 \bigg\}
	+ c_{10} \cdot \bigg\{ \io u^\frac{nq}{n+2q-2} dx \bigg\}^\frac{n+2q-2}{n}
  \eea
  for all $t\in (0,\tm)$,
  because $\sup_{t\in (0,\tm)} \io w(\cdot,t) dx<\infty$ by Lemma \ref{lem21} and (\ref{22.5}).\\
  Since
  \bas
	\frac{nq}{n+2q-2} -1 = \frac{(n-2)(q-1)}{n+2q-2}
  \eas
  clearly is positive, and since (\ref{22.87}) ensures that
  \bas
	\frac{\frac{nq}{n+2q-2}}{\frac{n(p+m-1)}{n-2}}
	= \frac{n}{n+2q-2} <1
  \eas
  and that
  \bas
	\frac{n(p+m-1)+2-n}{n\cdot \big(\frac{nq}{n+2q-2}-1\big)}
	= \frac{n+2q-2}{n} \cdot \frac{n(p+m-1)+2-n}{(n-2)(q-1)}
	= \frac{n+2q-2}{n},
  \eas
  we may once more draw on Lemma \ref{lem211} and (\ref{22.5}) to infer the existence of $c_{11}=c_{11}(p)>0$ satisfying
  \bas
	c_{10} \cdot \bigg\{ \io u^\frac{nq}{n+2q-2} dx \bigg\}^\frac{n+2q-2}{n}
	&=& c_{10} \cdot \bigg\{ \io u^\frac{nq}{n+2q-2} dx \bigg\}^\frac{n(p+m-1)+2-n}{n\cdot (\frac{nq}{n+2q-2}-1)} \\
	&\le& c_{11} \io \big|\na (u+1)^\frac{p+m-1}{2}\big|^2 dx
	+ c_{11}
	\qquad \mbox{for all } t\in (0,\tm).
  \eas
  Therefore, (\ref{22.9}) and (\ref{22.10}) show that
  \bas
	\frac{1}{q} \frac{d}{dt} \io w^q dx
	+ \lo \io w^q dx
	\le c_{11} \io \big|\na (u+1)^\frac{p+m-1}{2}\big|^2 dx
	+ \frac{4(q-1)}{q^2} + c_{11}
	\qquad \mbox{for all } t\in (0,\tm),
  \eas
  so that by (\ref{22.4}), (\ref{22.6}), (\ref{22.7}) and (\ref{22.8}),
  \bea{22.11}
	& & \hs{-20mm}
	\frac{d}{dt} \bigg\{ \io (u+1)^p dx
	+ \frac{c_1}{4q c_{11}} \io w^q dx \bigg\}
	+ \frac{c_1}{4} \io \big|\na (u+1)^\frac{p+m-1}{2}\big|^2 dx
	+ \frac{c_1 \lo}{4 c_{11}} \io w^q dx \nn\\
	&\le& c_8 \io w^\frac{p+m-1+\frac{2}{n}}{m+\frac{2}{n}-\sig} dx
	+ c_5 \io w^p dx
	+ c_{12}
	\qquad \mbox{for all } t\in (0,\tm)
  \eea
  with $c_{12}\equiv c_{12}(p):=\ko |\Om| + c_7 +  \frac{c_1}{4 c_{11}}\cdot \big(\frac{4(q-1)}{q^2} + c_{11}\big)$.
  Here we note that thanks to (\ref{22.3}) and the last restriction in (\ref{22.2}),
  \bas
	\frac{p+m-1+\frac{2}{n}}{m+\frac{2}{n}-\sig} < q
	\qquad \mbox{and} \qquad
	p=q-\Big(\frac{n(p+m-1)}{n-2} -p\Big)
	= q-\frac{2p-n(1-m)}{n-2} < q,
  \eas
  whence two further applications of Young's inequality provide $c_{13}=c_{13}(p)>0$ such that
  \be{22.12}
	c_8 \io w^\frac{p+m-1+\frac{2}{n}}{m+\frac{2}{n}-\sig} dx
	+ c_5 \io w^p dx
	\le \frac{c_1 \lo}{8c_{11}} \io w^q dx
	+ c_{13}
	\qquad \mbox{for all } t\in (0,\tm).
  \ee
  We finally observe that since, again due to the last inequality in (\ref{22.2}), we have $p  < \frac{n(p+m-1)}{n-2}$, we may once more
  combine Lemma \ref{lem211} with (\ref{22.5}) to see that with some $c_i=c_i(p)$, $i\in\{14,15\}$, we have
  \bas
	\frac{c_1}{4} \io \big|\na (u+1)^\frac{p+m-1}{2}\big|^2 dx
	\ge c_{14} \cdot \bigg\{ \io (u+1)^p dx \bigg\}^{\iota_1} - c_{15}
	\qquad \mbox{for all } t\in (0,\tm),
  \eas
  where $\iota_1:=\frac{n(p+m-1)+2-n}{n(p -1)}$ is positive due to the first condition in (\ref{22.2}).\abs
  Abbreviating $c_{16}\equiv c_{16}(p):=\frac{c_1}{4q c_{11}}$, from (\ref{22.11}) and (\ref{22.12}) we thus obtain that
  \be{22.13}
	y(t):=\io \big(u(\cdot,t)+1\big)^p dx
	+ c_{16} \io w^q(\cdot,t) dx,
	\qquad t\in (0,\tm),
  \ee
. satisfies
  \bas
	y'(t)
	+ c_{14} \cdot \bigg\{ \io (u+1)^p dx \bigg\}^{\iota_1}
	+ \frac{c_1 \lo}{8c_{11}} \io w^q dx
	\le c_{12} + c_{13} + c_{15}
	\qquad \mbox{for all } t\in (0,\tm),
  \eas
  so that writing $\iota_2:=\min\{\iota_1,1\}$ and estimating
  \bas
	y^{\iota_2}(t)
	&\le& 2^{\iota_2} \cdot \bigg\{ \io (u+1)^p dx \bigg\}^{\iota_2}
	+ (2c_{16})^{\iota_2} \cdot \bigg\{ \io w^q dx \bigg\}^{\iota_2} \\
	&\le& 2^{\iota_2} \cdot \Bigg\{ \bigg\{ \io (u+1)^p dx \bigg\}^{\iota_1} +1 \Bigg\}
	+ (2c_{16})^{\iota_2} \cdot \bigg\{ \io w^q dx + 1 \bigg\}
	\qquad \mbox{for all } t\in (0,\tm),
  \eas
  we arrive at the autonomous ODI
  \bas
	y'(t) + c_{17} y^{\iota_2}(t) \le c_{18}
	\qquad \mbox{for all } t\in (0,\tm)
  \eas
  with $c_{17}\equiv c_{17}(p):=\min \big\{ 2^{-\iota_2}  c_{14} \, , \, (2c_{16})^{-\iota_2}  \cdot \frac{c_1 \lo}{8c_{11}} \big\}$
  and $c_{18}\equiv c_{18}(p):= c_{12} + c_{13} + c_{15} + c_{17} \cdot \big\{ 2^{\iota_2} + (2c_{16})^{\iota_2} \big\}$.
  As thus
  \bas
	\io (u+1)^p dx
	\le \max \bigg\{ \io (u_0+1)^p dx \, , \, \Big(\frac{c_{18}}{c_{17}}\Big)^\frac{1}{\iota_2} \bigg\}
	\qquad \mbox{for all } t\in (0,\tm),
  \eas
  we infer that indeed (\ref{22.1})-(\ref{22.01}) holds when (\ref{kl2}) is satisfied.\abs
  In the case when (\ref{kl2}) fails, (\ref{22.1}) can be derived in quite a similar manner, with the only significant difference
  consisting in the way how the last summand in (\ref{22.4}) is dealt with.
  Then, namely, it is sufficient to simply estimate
  \bas
	p\kt \io (u+1)^{p-1} w dx
	\le (p-1) \kt \io (u+1)^p dx
	+ \kt \io w^p dx
	\qquad \mbox{for all } t\in (0,\tm),
  \eas
  and to use this as a substitute for (\ref{22.6}).
  By means of (\ref{mass_rough}), namely, we thus readily see that for each $T>0$ there exists $c_{19}=c_{19}(p,T)>0$
  such that for $y$ from (\ref{22.13}) we have
  \bas
	y'(t) \le c_{19}
	\qquad \mbox{for all } t\in (0,\tm),
  \eas
  and that hence we may conclude as intended also in this general case.
\qed
Our main result on global existence and boundedness in (\ref{0}) under the assumption (\ref{sig2}) has thereby essentially
been accomplished:\abs
\proofc of Proposition \ref{prop23}. \quad
  Based on a Moser-type iteration (cf. e.g. \cite[Lemma A.1]{taowin_subcrit}), from Lemma \ref{lem22}
  as well as standard parabolic and elliptic regularity theory applied to the third and second sub-problem in (\ref{0})
  it readily follows that for any $T>0$ there exists $c_1(T)>0$ such that
  \bas
	\|u(\cdot,t)\|_{L^\infty(\Om)} + \|v(\cdot,t)\|_{L^\infty(\Om)} + \|w(\cdot,t)\|_{L^\infty(\Om)} \le c_1(T)
	\qquad \mbox{for all } t\in (0,\tm) \cap (0,T),
  \eas
  with
  \bas
	\sup_{T>0} c_1(T)<\infty
	\qquad \mbox{if (\ref{kl2}) is satisfied.}
  \eas
  The claim thus results from Lemma \ref{lem_loc}.
\qed
\subsection{Criticality of the line $\sig=\frac{2}{n}$. Proof of Proposition \ref{prop25}}
In the framework of Proposition \ref{prop25}, drawing on dissipative effects of diffusion in the first equation from (\ref{0})
seems not expedient when $m$ attains very large negative values; in fact, it is well conceivable that unbounded solutions
exist in such cases.
That no such singularity formation can occur within finite time, however, will be a consequence of Lemma \ref{lem20}
when combined with a second application of Lemma \ref{lem21}, which provides the following major step toward Proposition \ref{prop25}.
\begin{lem}\label{lem24}
  Let $n\ge 3$ and $\Om\subset\R^n$ be a bounded domain with smooth boundary,
  and suppose that (\ref{kl}) holds, and that $D$ and $S$ satisfy (\ref{DSreg}), (\ref{Dl}) and (\ref{Su})
  with some $\xi_0>0, \kD>0, \KS>0, m\in\R$
  and $\sig\in\R$ such that (\ref{sig3}) is satisfied.
  Then assuming that $u_0$ and $w_0$ are such that (\ref{init}) holds and that in Lemma \ref{lem_loc} we have $\tm<\infty$,
  for any $p>1$ one can find $C(p)>0$ fulfilling
  \be{24.1}
	\|u(\cdot,t)\|_{L^p(\Om)} \le C(p)
	\qquad \mbox{for all } t\in (0,\tm).
  \ee
\end{lem}
\proof
 We evidently need to consider the case when $\sig\ge 0$ only, in which letting
  \be{24.02}
	k:=\frac{p}{p\sig+1-\sig}
	\qquad \mbox{and} \qquad
	l:=\frac{p}{1-\sig}
  \ee
  we obtain positive numbers $k$ and $l$ satisfying $l\ge 1$ and $k-1=\frac{(p-1)(1-\sig)}{p\sig+1-\sig}>0$ as well as
  \bas
	(n-2k)l-nk
	&=& \frac{np}{1-\sig} - \Big(\frac{2p}{1-\sig}+n\Big) \cdot \frac{p}{p\sig+1-\sig} \\
	&=& \frac{np(p\sig+1-\sig)-(2p+n-n\sig)p}{(1-\sig)(p\sig+1-\sig)} \\
	&=& \frac{(n\sig-2)p^2}{(1-\sig)(p\sig+1-\sig)} \\[2mm]
	&<& 0
  \eas
  according to our hypothesis $\sig<\frac{2}{n}$.
  We may therefore employ Lemma \ref{lem21} to find $c_1=c_1(p)>0$ such that
  \bas
	\|w\|_{L^\frac{p}{1-\sig}(\Om)} \le c_1 \cdot \bigg\{ 1 + \sup_{t'\in (0,t)} \|u(\cdot,t')\|_{L^k(\Om)} \bigg\}
	\qquad \mbox{for all } t\in (0,\tm),
  \eas
  so that since $k\in (1,p)$ and thus
  \bas
	\|u+1\|_{L^k(\Om)}
	\le \|u+1\|_{L^p(\Om)}^\frac{p(k-1)}{(p-1)k} \|u+1\|_{L^1(\Om)}^{1-\frac{p(k-1)}{(p-1)k}}
	\qquad \mbox{for all } t\in (0,\tm)
  \eas
  by the H\"older inequality, from (\ref{mass_rough}) we infer that as we are assuming that $\tm<\infty$, we can find $c_2=c_2(p)>0$
  such that
  \be{24.2}
	\|w\|_{L^\frac{p}{1-\sig}(\Om)} \le c_2 Y^\frac{k-1}{(p-1)k}(t)
	\qquad \mbox{for all } t\in (0,\tm),
  \ee
  where
  \bas
	Y(t):=\max_{t'\in [0,t]} y(t'),
	\qquad t\in [0,\tm),
  \eas
  with
  \bas
	y(t):=1+\io \big(u(\cdot,t)+1\big)^p dx,
	\qquad t\in [0,\tm).
  \eas
  We now combine Lemma \ref{lem20} with Young's inequality to see that with some $c_3=c_3(p)>0$ and $c_4=c_4(p)>0$ we have
  \bea{24.3}
	y'(t)
	&\le& c_3 \io (u+1)^{p+\sig-1} w dx + c_3 \io (u+1)^{p-1} w dx \nn\\
	&\le& 2c_3 \io (u+1)^{p+\sig-1} w dx + c_4
	\qquad \mbox{for all } t\in (0,\tm),
  \eea
  because again due to Lemma \ref{lem21} and (\ref{mass_rough}), our assumption on finiteness of $\tm$ clearly implies that
  $\sup_{t\in (0,\tm)} \io w(\cdot,t) dx <\infty$.
  But by the H\"older inequality and (\ref{24.2}),
  \bas
	\io (u+1)^{p+\sig-1} w dx
	&\le& \bigg\{ \io (u+1)^p dx \bigg\}^\frac{p+\sig-1}{p} \cdot \|w\|_{L^\frac{p}{1-\sig}(\Om)} \\
	&\le& c_2 y^\frac{p+\sig-1}{p}(t) \cdot Y^\frac{k-1}{(p-1)k}(t) \\
	&\le& c_2Y^{\frac{p+\sig-1}{p} + \frac{k-1}{(p-1)k}}(t)
	\qquad \mbox{for all } t\in (0,\tm),
  \eas
  so that from (\ref{24.3}) and the inequality $Y\ge 1$ we infer that
  \bas
	y'(t) \le c_5 Y(t)
	\qquad \mbox{for all } t\in (0,\tm)
  \eas
  with $c_5\equiv c_5(p):=2c_2 c_3 + c_4$,
  because (\ref{24.02}) entails that
  \bas
	\frac{p+\sig-1}{p} + \frac{k-1}{(p-1)k}
	&=& 1 + \frac{\sig-1}{p} + \frac{1}{p-1} \cdot \Big(1-\frac{p\sig+1-\sig}{p}\Big) \\
	&=& 1 + \frac{\sig-1}{p} + \frac{p-p\sig-1+\sig}{p(p-1)} \\
	&=& 1.
  \eas
  An integration now shows that
  \bas
	y(t)
	&\le& y(0) + c_5 \int_0^t Y(\tau) d\tau
	\qquad \mbox{for all } t\in (0,\tm)
  \eas
  and that thus, by definition of $Y$,
  \bas
	Y(t) \le y(0) + c_5 \int_0^t Y(\tau) d\tau
	\qquad \mbox{for all } t\in (0,\tm).
  \eas
  Relying on the continuity of $Y$, we may therefore draw on Gronwall's inequality to conclude that
  \bas
	Y(t) \le y(0) e^{c_5 t} \le y(0) e^{c_5 \tm}
	\qquad \mbox{for all } t\in (0,\tm)
  \eas
  which implies (\ref{24.1}).
\qed
Indeed, this immediately implies global solvability whenever (\ref{sig3}) holds:\abs
\proofc of Proposition \ref{prop25}. \quad
  In contrast to the argument in Lemma \ref{lem24},
  we now make essential use of the presence of diffusion and the algebraic lower estimate in (\ref{Dl}).
  Together with the outcome of Lemma \ref{lem24}, namely, again thanks to parabolic and elliptic regularity
  arguments this enables us to conclude from a Moser-type reasoning (\cite[Lemma A.1]{taowin_subcrit}) that
  if $\tm$ was finite, then we could find $c_1>0$ such that
  $\|u(\cdot,t)\|_{L^\infty(\Om)} \le c_1$ for all $t\in (0,\tm)$.
  As this would contradict the extensibility criterion (\ref{ext}) in Lemma \ref{lem_loc}, however, this
  already completes the proof.
\qed

\bigskip

{\bf Acknowledgement.} \quad
 Youshan Tao was supported by the {\em National Natural Science Foundation of China
   (No. 12171316)}.  Michael Winkler acknowledges support of the
  {\em Deutsche Forschungsgemeinschaft} (Project No.~462888149).
\end{document}